\documentclass[11pt]{amsart}

\usepackage[utf8]{inputenc} 	
\usepackage[T1]{fontenc}
\usepackage{lmodern}
\usepackage{yhmath}

\usepackage{amsthm,amsmath,amssymb,amsbsy,bbm,mathrsfs,supertabular,eurosym,graphicx,enumitem,xcolor}
\usepackage{mathtools}
\usepackage{leftidx} 
\usepackage{scalerel}

\newtheoremstyle{BBstyle0}  {}{}{\itshape}{}{\bfseries}{}{6pt}{}
\newtheoremstyle{BBstyle1}  {3pt}{3pt}{\rmfamily}{}{\itshape}{: }{3pt}{}
\newtheoremstyle{BBstyle2}  {3pt}{3pt}{\itshape}{}{\bfseries\large}{}{0pt}{}
\newtheoremstyle{BBstyle3}  {}{}{\itshape}{}{\bfseries}{: }{3pt}{}
\newtheoremstyle{BBstyle4}  {}{}{\rmfamily}{}{\bfseries}{}{6pt}{}
  
\usepackage[normalem]{ulem} 

\usepackage[authoryear]{natbib}
\usepackage{enumitem}
\usepackage[hidelinks]{hyperref}  

\newcommand{\pa}[1]{\left({#1}\right)}
\newcommand{\norm}[1]{\left\|{#1}\right\|}
\newcommand{\cro}[1]{\left[{#1}\right]}
\newcommand{\ab}[1]{\left|{#1}\right|}
\newcommand{\ac}[1]{\left\{{#1}\right\}}

\newcommand{\argmin}{\mathop{\rm argmin}}

\newcommand{\Var}{\mathop{\rm Var}\nolimits}

\newcommand{\E}{{\mathbb{E}}}
 
\renewcommand{\L}{{\mathbb{L}}}
\newcommand{\N}{{\mathbb{N}}}
\renewcommand{\P}{{\mathbb{P}}}
 
\newcommand{\R}{{\mathbb{R}}}
\renewcommand{\S}{{\mathbb{S}}} 
\newcommand{\T}{{\mathbb{T}}} 
\newcommand{\Z}{{\mathbb{Z}}}
\renewcommand{\H}{{\mathbb{H}}}

\newcommand{\U}{{\mathbb{U}}}

\newcommand{\sA}{{\mathscr{A}}}
\newcommand{\sB}{{\mathscr{B}}}
\newcommand{\sC}{{\mathscr{C}}}
\newcommand{\sD}{{\mathscr{D}}}

\newcommand{\sP}{{\mathscr{P}}}

\newcommand{\fra}{{\mathfrak{a}}}

\newcommand{\frc}{{\mathfrak{c}}}
\newcommand{\frd}{{\mathfrak{d}}}

\newcommand{\frm}{{\mathfrak{m}}}

\newcommand{\frs}{{\mathfrak{s}}} 

\newcommand{\fru}{{\mathfrak{u}}}

\newcommand{\frJ}{{\mathfrak{J}}}

\newcommand{\frM}{{\mathfrak{M}}}

\DeclareMathAlphabet{\mathscrbf}{OMS}{mdugm}{b}{n}
\newcommand{\cA}{{\mathcal{A}}}
\newcommand{\cB}{{\mathcal{B}}}

\newcommand{\cE}{{\mathcal{E}}}

\newcommand{\cG}{{\mathcal{G}}}

\newcommand{\cJ}{{\mathcal{J}}}
\newcommand{\cK}{{\mathcal{K}}}
 
\newcommand{\cM}{{\mathcal{M}}}
\newcommand{\cN}{{\mathcal{N}}}
\newcommand{\cO}{{\mathcal{O}}}

\newcommand{\cS}{{\mathcal{S}}} 
\newcommand{\cT}{{\mathcal{T}}}
\newcommand{\cU}{{\mathcal{U}}}

\newcommand{\cX}{{\mathcal{X}}}

\newcommand{\gf}{{\mathbf{f}}}
\renewcommand{\gg}{{\mathbf{g}}} 
\newcommand{\gh}{{\mathbf{h}}}

\newcommand{\gk}{{\mathbf{k}}}
\newcommand{\gl}{{\mathbf{l}}}

\newcommand{\gr}{{\mathbf{r}}}
\newcommand{\gs}{{\mathbf{s}}} 
\newcommand{\gt}{{\mathbf{t}}}
\newcommand{\gu}{{\mathbf{u}}}
\newcommand{\gv}{{\mathbf{v}}}
\newcommand{\gw}{{\mathbf{w}}}
\newcommand{\gx}{{\mathbf{x}}}

\newcommand{\gS}{{\mathbf{S}}} 
\newcommand{\gT}{{\mathbf{T}}}

\newcommand{\bs}[1]{\boldsymbol{#1}}

\newcommand{\bsX}{{\bs{X}}}

\newcommand{\gmu}{{\bs{\mu}}}

\newcommand{\gpi}{{\bs{\pi}}}

\newlist{lista}{enumerate}{2}
\setlist[lista,1]{label=\alph*),ref=\alph*)}

\newlist{listi}{enumerate}{1}
\setlist[listi,1]{label=(\roman*),ref=(\roman*),align=left}

\newlist{lists}{enumerate}{3}
\setlist[lists,1]{label=$\blacktriangleright$,align=right}
\setlist[lists,2]{label=$\bullet$,align=right}
\setlist[lists,3]{label=$\star$,align=right}

\newcommand{\eref}[1]{(\ref{#1})}

\renewcommand{\ge}{\geqslant}
\renewcommand{\le}{\leqslant}
\newcommand{\1}{1\hskip-2.6pt{\rm l}}
\newcommand{\0}{{\bf 0}}

\newcommand{\scal}[2]{\left\langle #1,#2\right\rangle}

\newcommand{\etc}[1]{#1_1,\ldots,#1_n}

\newcommand{\on}{^{\otimes n}}
\newcommand{\tr}{^{\top}}
\newcommand{\et}{^{\star}}
\def\vide{{\varnothing}}

\newcommand{\map}[5]{
\begin{array}{l|rcl}
#1: & #2 & \longrightarrow & #3 \\
    & #4 & \longmapsto & #5 
\end{array}
    }

\DeclarePairedDelimiter\floor{\lfloor}{\rfloor}

\newtheorem{thm}{Theorem}
\newtheorem{lem}{Lemma}
\newtheorem{prop}{Proposition}
\newtheorem{df}{Definition}
\newtheorem{cor}{Corollary}

\theoremstyle{definition}
\newtheorem{exa}{Example}

\usepackage[english]{babel}

\def\gell{{\bs{\ell}}}

\def\ray{{r}}

\def\X{\mathbb{X}}

\begin{document}
\title{Robust Bayes-like estimation of the distribution of a determinantal process}
\author{Yannick BARAUD}
\address{\parbox{\linewidth}{Department of Mathematics,\\
University of Luxembourg\\
Maison du nombre\\
6 avenue de la Fonte\\
L-4364 Esch-sur-Alzette\\
Grand Duchy of Luxembourg}}
\email{yannick.baraud@uni.lu}
\date{November 6, 2026}

\begin{abstract}
We consider the problem of estimating the distribution of a determinantal process (DP) on a
measured space from the observation of $n$ independent copies of it. Our estimator is obtained by
sampling from a $T$-posterior distribution associated with a suitable prior, and we establish
non-asymptotic risk bounds for the Hellinger loss which hold with high probability. These bounds
remain valid when the data are only independent rather than i.i.d., and when the model is
misspecified, which makes the estimator robust, for instance, to the presence of a small fraction
of outliers. For determinantal projection processes of rank $k$ whose kernels are associated with
$k$-dimensional subspaces of a $D$-dimensional space, Haar priors on Grassmannians lead to a
rate of order $\sqrt{k(D-k)\log n/n}$, and model selection extends this result to infinite-dimensional
settings. For general DPs, we introduce priors based on the spectral decomposition of the kernel
that combine uniform distributions on the eigenvalues with Haar measures on Stiefel manifolds.
When the underlying space is finite with $N$ elements, we obtain a rate of order $N\sqrt{\log n/n}$,
without any assumption on the kernel, and derive that the $2^{N}$ principal minors of the kernel can
be estimated simultaneously. Finally, for translation-invariant DPs on the torus $\mathbb{T}^d$,
we obtain rates that only depend on the smoothness of the kernel at the origin and that
are adaptive with respect to it. Combined with a new lower bound on the Hellinger distance between such
DPs, these rates also yield estimators of the kernel in $\mathbb{L}_2$.
\end{abstract}
\maketitle

\section{Introduction}

We observe $\bsX=(X_{1},\ldots,X_{n})$, where the $X_{i}$ are $n$ finite point processes on a measure space $(\cX,\cA,\mu)$, that is, for every $i\in\{1,\ldots,n\}$, $X_{i}$ takes its values in the class $\X$ of all finite subsets of $\cX$. Our aim is to estimate the distribution of these data, assuming that it is close, in a sense to be specified later, to that of an $n$-sample drawn from a determinantal distribution on $\X$.

\subsection{Determinantal processes}
We equip $\X$ with a $\sigma$-algebra $\Xi$ that makes all the mappings
\[
\map{i_{A}}{\X}{\N}{\gx}{\ab{\gx\cap A}={\rm Card}(\gx\cap A)}
\]
measurable for every $A\in\cA$. In particular, the set $\X_{k}=i_{\cX}^{-1}(\{k\})$, which is the class of all subsets of $\cX$ with cardinality $k$, is measurable for every $k\in\N$. The measurable space $(\X,\Xi)$ is endowed with the measure $\gmu$ defined by 
\begin{align*}
\int_{\X}f(x)d\gmu(x)=f(\vide)+\sum_{k\ge 1}\frac{1}{k!}\int_{\cX^{k}}f(\{x_{1},\ldots,x_{k}\})\1_{\X_{k}}\pa{\{x_{1},\ldots,x_{k}\}}d\mu(x_{1})\ldots d\mu(x_{k}),
\end{align*}
for every measurable function $f\colon\X=\bigcup_{k\ge 0}\X_{k}\to \R_{+}$.

The determinantal point processes (DPP for short) that we consider in this paper are (possibly infinite) mixtures of {\em determinantal projection processes} (PDPs for short), which form a subclass of determinantal processes. We first define the latter.
Given a linear subspace $U$ of $\H=\L_{2}(\cX,\cA,\mu)$ of finite dimension $k\ge 1$ and an orthonormal basis $\gu=\{u_{1},\ldots,u_{k}\}$ of $U$, we define the function $p_{U}$ which associates with every $x=\{x_{1},\ldots,x_{k}\}\in \X_{k}$ the value
\begin{equation}\label{def-dens}
p_{U}(x)=\det\pa{\sum_{i=1}^{k}u_{j}(x_{i}) u_{j'}(x_{i})}_{1\le j,j'\le k}
=\ab{\det\left(
\begin{array}{ccc}
u_{1}(x_{1}) && u_{k}(x_{1}) \\
\vdots & \vdots & \vdots \\
u_{1}(x_{k}) &  & u_{k}(x_{k})
\end{array}
\right)}^{2}.
\end{equation}
It is easy to check that $p_{U}$ depends only on the linear space $U$ and not on the choice of an orthonormal basis $\gu$. In fact, $p_{U}$ is the density, with respect to $\gmu$, of a probability distribution on $\X_{k}$, which we denote by $P_{U}$. A random variable with such a density is called {\em a determinantal projection process} or, more precisely, {\em a $k$-determinantal projection process}. Note that the observation $X$ of a $k$-PDP automatically determines the value of $k$: since $X$ belongs to $\X_{k}$, its cardinality $|X|$ is necessarily $k$.  

If the functions $u_{1},\ldots,u_{k}$ are continuous and the points $x_{1},\ldots,x_{k}$ are close enough to one another, the value of the determinant in \eref{def-dens} is small, and so is the density.  Projection determinantal processes (as well as general DPPs) thus generate point patterns that exhibit repulsion.

A simple example of a $k$-PDP is the following. If $U$ is the linear span of $k$ functions $u_{1},\ldots,u_{k}$  lying in the unit sphere $\S$ of $\H$ with disjoint supports $A_{1},\ldots,A_{k}$ respectively, then for $x=\{x_{1},\ldots,x_{k}\}$ with $x_{i}\in A_{i}$ for $i\in\{1,\ldots,k\}$, 
\[
p_{U}(x)=u_{1}^{2}(x_{1})\ldots u_{k}^{2}(x_{k}).
\]
The distribution of $X$ is then the same as that of $\{X_{1},\ldots,X_{k}\}$ where $X_{1},\ldots,X_{k}$ are independent random variables with densities $u_{1}^{2},\ldots,u_{k}^{2}$ respectively. In particular, the problem of estimating the distribution of $X$ is equivalent to that of estimating the densities $u_{1}^{2},\ldots,u_{k}^{2}$ simultaneously. This example helps  relate our bounds for estimating the density $p_{U}$ to more classical ones for estimating a density of the form $p=u^{2}$ under the assumption that $u$ belongs to a linear space.  

In the literature, the distribution of a general DPP is associated with a self-adjoint nonnegative operator on $\H$ of the form
\begin{equation}\label{eq-oper}
f\mapsto \int_{\cX}K(x,y)f(y)d\mu(y)
\end{equation}
where $K$ denotes a suitable kernel on $\cX\times\cX$. Hereafter, we only consider symmetric kernels $K$, that is, $K(x,y)=K(y,x)$ for every $x,y\in\cX$. Following the probabilistic construction of DPPs described in Hough {\em et al.}~\citeyearpar{MR2216966}, we may consider only kernels $K$ of the form
\begin{equation}\label{def-K}
K(x,y)=\sum_{j\ge 1}l_{j}u_{j}(x)u_{j}(y)\quad \text{for every $x,y\in \cX$}
\end{equation}
where the $l_{j}$ belong to $[0,1]$ and satisfy $\sum_{j\ge 1}l_{j}<+\infty$, and $(u_{j})_{j\ge 1}$ is an orthonormal system of functions in $\H$. Under these assumptions, the operator \eref{eq-oper} is well defined and 
\begin{equation}\label{eq-tracefinie}
\int_{\cX}K(x,x)d\mu(x)=\sum_{j\ge 1}l_{j}<+\infty. 
\end{equation}
According to this construction, a determinantal process $X$ with distribution $P_{K}$ associated with $K$ can be built as follows. First generate a sequence of independent Bernoulli random variables $\epsilon_{j}$ with parameters $l_{j}$ for every $j\ge 1$. Since $\sum_{j\ge 1}l_{j}<+\infty$, the Borel--Cantelli lemma implies that the set $J=\{j\ge 1,\ \epsilon_{j}=1\}$ is almost surely finite, and so is the dimension of the linear subspace $U_{J}\subset \H$ generated by the vectors $\{u_{j}, j\in J\}$. Finally, conditionally on $J$, draw $X$ from the PDP distribution $P_{U_{J}}$ associated with $U_{J}$ when $J\ne \vide$, and set $X=\vide$ otherwise. In the case of a PDP, $K$ is given by $K(x,y)=\sum_{j=1}^{k}u_{j}(x)u_{j}(y)$ and the operator \eref{eq-oper} is simply the orthogonal projector onto the linear space $U$ in $\H$. This construction shows that $P_{K}$ is a (possibly countable) mixture of determinantal projection processes, and 
\begin{equation}\label{eq-nb-pt}
\E\cro{\ab{X}}=\int_{\cX}K(x,x)d\mu(x).
\end{equation}
\subsection{Our estimation strategy and contributions}
Our framework fits into the more general one of estimating the distribution of $n$ observed random variables $\bsX=(\etc{X})$  
by means of a statistical model of the form $\{P_{\theta}\on,\; \theta\in \Theta\}$ where $\sP=\{P_{\theta}=p_{\theta}\cdot\gmu\}$ is a dominated family of probabilities on $(\X,\Xi,\gmu)$. Note that this statistical model assumes that the data are i.i.d.\ with a distribution in $\sP$, although this may not be the case in reality. The problem of estimating the distribution of the data based on a possibly misspecified model has been tackled in Baraud~\citeyearpar{baraud2021robust} and Baraud and Romon~\citeyearpar{Baraud:2026aa}. We use here the theory of $T$-Bayes estimation developed there to address our estimation problem. 

We present this theory in the simple case where the data $X_{1},\ldots,X_{n}$ are independent, although Baraud and Romon~\citeyearpar{Baraud:2026aa} extended it to data with some form of dependence. Let us therefore assume that $X_{1},\ldots,X_{n}$ are independent (but not necessarily i.i.d.) with marginal distributions $P_{1},\ldots,P_{n}$ respectively. If our statistical model $\{P_{\theta}\on,\; \theta\in \Theta\}$ provides a reasonable approximation of reality, there exists a distribution $P_{\theta}$ in $\sP$ which is close to these marginal distributions, or at least to most of them. To measure this proximity (and, later on, the performance of our estimator), we use the Hellinger distance. We recall that the Hellinger distance between two probabilities $P$ and $Q$ dominated by $\nu$ on a space $E$ is defined by 
\[
h\pa{P,Q}=\sqrt{\frac{1}{2}\int_{E}\pa{\sqrt{\frac{dP}{d\nu}}-\sqrt{\frac{dQ}{d\nu}}}^{2}d\nu},
\]
the result being independent of the choice of $\nu$. It is proven in Baraud~\citeyearpar{baraud2021robust} and Baraud and Romon~\citeyearpar{Baraud:2026aa} that, given a prior $\gpi$ on the parameter space $\Theta$, it is possible to design a posterior distribution $\gpi_{\bsX}$ based on $\bsX$ which satisfies the following property. For every $\xi>0$ and $\theta\in \Theta$, with probability at least $1-2e^{-\xi}$, a (randomised) estimator $\widehat \theta=\widehat \theta(\bsX)$ drawn from $\gpi_{\bsX}$ satisfies 
\begin{equation}\label{eq-riskbase}
h^{2}\pa{P_{\theta},P_{\smash{\widehat \theta}}}\le C\cro{h^{2}\pa{\frac{1}{n}\sum_{i=1}^{n}P_{i},P_{\theta}}+\frac{r(\gpi,\theta)+1+\xi}{n}}
\end{equation}
where $C$ is a positive constant and $r(\gpi,\theta)$ a dimensional term that depends on the behaviour of the prior $\gpi$ in a neighbourhood of $\theta$. We call $r(\gpi,\theta)$ {\em the $\gpi$-complexity at $\theta$}, $\gpi_{\bsX}$ {\em the $T$-Bayes posterior} and $\widehat \theta=\widehat \theta(\bsX)$ a $T$-Bayes estimator. An inequality of the form~\eref{eq-riskbase} cannot in general be achieved by means of a classical Bayesian approach, and we refer the reader to Baraud and Birg\'e~\citeyearpar{BarBir2020} for a discussion of the instability of the Bayesian procedure under misspecification. 

Interestingly, the $\gpi$-complexity function $\theta\mapsto r(\gpi,\theta)$ depends only on the prior $\gpi$ and the statistical model, not on the true distribution of the data. It is therefore a known quantity. In the special situation where the parameter space $\Theta$ is a convex compact subset of $\R^{D}$ and $\gpi$ is the uniform distribution on $\Theta$, $r(\gpi,\theta)$ is typically not larger than $D$, up to a constant (see Baraud and Birg\'e~\citeyearpar{BarBir2020} and Proposition~1 in  Baraud and Romon~\citeyearpar{Baraud:2026aa} for an extension). These priors have the advantage of being agnostic with respect to the location of the parameter within the parameter space. 

Integrating \eref{eq-riskbase} with respect to $\xi$ yields the risk bound 
\begin{equation}\label{eq-risk00}
\E\cro{h^{2}\pa{P_{\theta},P_{\smash{\widehat \theta}}}}\le C\cro{h^{2}\pa{\frac{1}{n}\sum_{i=1}^{n}P_{i},P_{\theta}}+\frac{r(\gpi,\theta)+1}{n}}.
\end{equation}
It is of order $r(\gpi,\theta)/n$ when the data are truly i.i.d.\ with distribution $P_{\theta}$, that is, when the statistical model is exact, while, in the case of misspecification, it is only inflated by the extra term $h^{2}\pa{n^{-1}\sum_{i=1}^{n}P_{i},P_{\theta}}$. Since 
\[
h^{2}\pa{n^{-1}\sum_{i=1}^{n}P_{i},P_{\theta}}\le \frac{1}{n}\sum_{i=1}^{n}h^{2}(P_{i},P_{\theta})
\]
we see that this additional term remains small compared to $r(\gpi,\theta)/n$ if most of the marginal distributions are close enough to $P_{\theta}$. In a parametric setting where $\Theta$ is a subset of $\R^{D}$ and $r(\gpi,\theta)$ is of order $D$, the risk bound therefore remains of the same order as for i.i.d.\ data with distribution $P_{\theta}\in\sP$, even if the distributions of $n'$ of these data were replaced by arbitrary ones, as long as $n'$ remains small compared to $D$. 

Since the risk of our $T$-Bayes estimators is only inflated by an approximation term in the case of misspecification, it is enough to study them when the model is well specified. Throughout this paper, we  therefore restrict ourselves to the situation where $\etc{X}$ are i.i.d.\ with a determinantal distribution $P_{\theta\et}$. 

In the case of a $k$-PDP, the parameter $\theta\et$ takes the form of a linear subspace $U$ of $\H$ of dimension $k$. The parameter space is therefore the Grassmannian $\cG_{k}(\H)$ (extended to Hilbert spaces)  of all $k$-dimensional subspaces of $\H$. This unusual parameter space raises the following questions: 
\begin{itemize}
\item Is there a geometric distance on this Grassmannian that makes it possible to bound the Hellinger distance between $P_{\theta}$ and $P_{\theta'}$ by the distance between the corresponding linear spaces? 
\item Are there agnostic priors on this Grassmannian that play a role similar to that of the uniform priors on compact subsets of $\R^{D}$? 
\end{itemize}
As for the first question, we show that a distance based on the principal angles between two $k$-dimensional linear spaces $\theta$ and $\theta'$ actually controls the Hellinger distance between $P_{\theta}$ and $P_{\theta'}$. As for the second one, Haar measures on some specific Grassmannian submanifolds provide ``flat'' priors that can be regarded as agnostic, since they are invariant under the action of the orthogonal group.    

As for more general DPPs, they can be parameterised either by kernels $K$ satisfying \eref{def-K} or, alternatively, by the parameter $\theta=((l_{j})_{j\ge 1}, (u_{j})_{j\ge 1})$ corresponding to the sequences of eigenvalues and eigenvectors of the 
operator \eref{eq-oper}. We then propose priors on the space of sequences as well as on Stiefel manifolds of orthonormal systems in $\H$, in order to control the complexity term $r(\gpi,\theta)$. 

A common feature of our approach is that we do not assume that the true parameter $\theta\et$ belongs to the support ${\rm supp}(\gpi)$ of the priors $\gpi$ that we propose. In particular, the classical Bayesian theory under misspecification developed in Kleijn and van der Vaart~\citeyearpar{MR2283395}, which describes the concentration of the posterior around the Kullback--Leibler projection of the true distribution onto the support of the prior, may not apply in this context: the true distribution $P_{\theta\et}$ is not in general be dominated by any element of $\{P_{\theta},\; \theta\in {\rm supp}(\gpi)\}$, so that the Kullback--Leibler divergence may be infinite. We rather use priors supported on subsets of the parameter set $\Theta$ and use the following fundamental property, which is a consequence of \eref{eq-risk00} and the triangle inequality, 
\begin{equation}\label{eq-risk01}
\E\cro{h^{2}\pa{P_{\theta\et},P_{\smash{\widehat \theta}}}}\le C\inf_{\theta\in {\rm supp}(\gpi)}\cro{h^{2}\pa{P_{\theta\et},P_{\theta}}+\frac{r(\gpi,\theta)+1}{n}}.
\end{equation}
We may therefore make our risk bound small by choosing a prior $\gpi$ for which the submodel $\{P_{\theta},\; \theta\in {\rm supp}(\gpi)\}\subset \sP$ provides a good enough approximation of $P_{\theta\et}$ (in Hellinger distance) while keeping the $\gpi$-complexity of ${\rm supp}(\gpi)$ to reasonable level. Throughout this paper, we propose such submodels, describe their approximation properties and evaluate the $\gpi$-complexities of the corresponding parameter spaces. This leads to various risk bounds. Some of them are optimal up to a logarithmic factor, as they generalise well-known bounds that have been obtained in more classical density estimation contexts.  

Our general results also apply to the case where $\cX=\{1,\ldots,N\}$ is finite, in which case $K$ can be identified with an $N\times N$ positive semi-definite matrix with eigenvalues in $[0,1]$. With our approach, we show that we can build an estimator $\widehat K$ of $K$ that satisfies some some numerical constant $C>0$, 
\[
\sup_{K}\E\cro{\sup_{J\subset \cX,J\ne \vide}\ab{\det K_{J}-\det \widehat K_{J}}}\le C N \sqrt{\frac{\log(n+1)}{n}}
\]
where the supremum is taken over positive semi-definite matrices with eigenvalues in $[0,1]$. In particular, $\widehat K$ can estimate the $2^{N}$ principal minors of $K$ simultaneously at the rate $1/\sqrt{n}$ up to a logarithmic factor. 

Another application concerns translation-invariant kernels on the $d$-dimensional torus $\T^{d}$. We assume that $K$ is of the form $K(x,y)=\overline K(x-y)$ for $x,y\in \T^{d}$, where $\overline K$ is an unknown even function on $\T^{d}$. In this framework, we show that the parameter $\overline K^{2}$ is identifiable. In particular, the function $\overline K$ is also identifiable when it is nonnegative on $\T^{d}$. Our approach allows us to establish non-asymptotic risk bounds with respect to the $\L_{2}$-loss for estimating $\overline K^{2}$. Somewhat surprisingly, these bounds depend only on the smoothness of $\overline K$ at the origin. As a corollary, we show that when $\overline K$ is bounded from below by some positive number $\fra >0$ on $\T^{d}$, our approach leads to an estimator $\widehat K$ of $\overline K$ that satisfies 
\begin{align*}
\E\cro{\norm{\overline K-\widehat K}_{2}^{2}}\le C'\cro{\frac{\overline K(0)\pa{1\vee \overline K(0)}}{\fra}}^{2}\inf_{D\ge 1}\cro{w\pa{\frac{1}{\pi D}}+\frac{(1+2D)^{d}}{n}\log n}
\end{align*}
where $w$ denotes the modulus of continuity of $\overline K$ at 0 and $\norm{\cdot}_{2}$ the $\L_{2}$-norm. In particular, if $\overline K$ is $(L,\alpha)$-H\"older continuous at 0 with $L,\alpha>0$, then $\widehat K$ estimates $\overline K$ at the rate $L^{d/(\alpha+d)}\pa{(\log n)/n}^{\alpha/(\alpha+d)}$.

\subsection{The literature}
From a theoretical point of view, few results are available on the estimation of the distribution of a determinantal process. When $\cX$ is a finite set, the problem of identifying the positive semi-definite matrix $K$ can only be solved up to $D$-similarity, that is, $K$ and $K'$ are associated with the same DPP distribution if and only if $K'=DKD$ where $D$ is a diagonal matrix with entries in $\{\pm 1\}$.
Rising {\em et al.}~\citeyearpar{RisingJustin2015Aeaf} (see also Kulesza and Taskar~\citeyearpar{10.1561/2200000044}) provided an efficient algorithm for recovering a symmetric matrix from its principal minors. In this framework, the Maximum Likelihood Estimator (MLE) of $L=K(I-K)^{-1}$ was studied in Brunel {\em et al.}~\citeyearpar{Brunel:2017aa} (see also Brunel {\em et al.}~\citeyearpar{Deter-rate-Brunel}). They showed that the MLE is asymptotically normal at the parametric rate when $L$  is irreducible, whereas its rate of convergence may be slower than the usual $1/\sqrt{n}$ rate otherwise. Urschel {\em et al.}~\citeyearpar{Urschel:2017aa} established the optimal sample complexity for estimating the matrix $K$, which is governed by the cycle sparsity of the graph associated with $K$, and proposed a fast algorithm based on the method of moments. Their results rely on the assumption that the non-zero entries of $K$ are bounded away from 0 by a known constant. 

For parametric stationary DPPs on $\R^{d}$, Biscio and Lavancier~\citeyearpar{MR3619702} studied minimum contrast estimators, based on Ripley's $K$-function or on the pair correlation function, of the type considered in Lavancier {\em et al.}~\citeyearpar{MR3382600}. A single DPP is observed on a bounded domain $D_{n}$ that grows to $\R^{d}$ with $n$, and they established the consistency and the asymptotic normality of these estimators as $n$ tends to infinity. Estimating-function methods based on the joint intensities of the process, which also cover the non-stationary case, are studied in Lavancier {\em et al.}~\citeyearpar{Lavancier_2020}. 

These results address the problem of parameter estimation in a situation where the model is well specified. With the exception of the work of Urschel {\em et al.}~\citeyearpar{Urschel:2017aa}, all these results are asymptotic. None of them considers the problem of estimating the distribution of the process. The only exception we are aware of is Baraud~\citeyearpar{Baraud:2012qy}. He establishes non-asymptotic risk bounds for the Hellinger loss. He uses an aggregation procedure based on robust tests, in the spirit of the $T$-estimators introduced by Birg\'e (see Birg\'e~\citeyearpar{MR2219712}). His approach relies  on a suitable discretisation of the parameter space. 

To our knowledge, estimating the parameters or the distribution of a DPP from a theoretical point of view by means of a Bayesian-type approach has not been addressed in the literature so far. This is the purpose of the present paper.

\subsection{Organisation of the paper and notation}
The paper is organised as follows. We describe our Bayes-like estimation strategy as well as its performance in Section~\ref{sect-strat}. We apply it to estimate the distribution of a PDP in Section~\ref{sect-PDP}.  Section~\ref{sect-DPP}
deals with general DPPs, including the case of a finite set $\cX$, and Section~\ref{sect-tore} with translation-
invariant DPPs on $\T^{d}$. In each of these sections, we relate the Hellinger distance between the distributions to a distance between the parameters and analyse the $\gpi$-complexity of some priors, which leads to various risk bounds. Finally, Section~\ref{sect-proof} is devoted to the proofs of our results.

Throughout this paper, it is convenient to identify a finite subset $\gx$  of $\cX$ with the corresponding counting measure, that is, the  measure that puts a mass of one at each element of  $\gx$. In particular, for a subset $A$ of $\cX$, the cardinality $|\gx\cap A|$ of the finite set $\gx\cap A$ is also equal to $\gx(A)$ when $\gx$ is viewed as a discrete measure. Given $t,s\in\R$, $s\wedge t$ and $s\vee t$ denote the minimum and the maximum of $s$ and $t$, respectively. In particular, $\log_{+}(t)=(\log t)\vee 1$ for $t>0$ and $(t)_{+}=t\vee 0$ for $t\in\R$. We set $\N\et=\N\setminus\{0\}$ and, for $x\in\R^{d}$, we denote by $|x|$ its Euclidean norm, $|x|_{1}=\sum_{i=1}^{d}|x_{i}|$ its $\ell_{1}$-norm and by $|x|_{\infty}=\max_{i=1,\ldots,d}|x_{i}|$ its supremum norm. Finally, $C,C',\ldots$ denote positive numerical constants whose values may vary from line to line. 

\section{Our estimation strategy and its performance}\label{sect-strat}
\subsection{$T$-Bayes posterior and $T$-Bayes estimator}\label{sect-TPost}
As already mentioned, we use a $T$-Bayes approach as described in Baraud~\citeyearpar{baraud2021robust} and Baraud and Romon~\citeyearpar{Baraud:2026aa} to address our estimation problem. Let $\bsX=(\etc{X})$ be an $n$-tuple of independent random variables with values in $\X$, $\sP=\{P_{\theta}=p_{\theta}\cdot\gmu, \theta\in \Theta\}$ a dominated family of distributions on $\X$ and a prior $\gpi$ on the parameter space $\Theta$. We define a posterior distribution $\gpi_{\!\bsX}$, based on the data $\bsX$, as follows. For $\theta,\theta'\in \Theta$ we set 
\[
\gT(\bsX,\theta,\theta')=\sum_{i=1}^{n}\psi\pa{\sqrt{\frac{p_{\theta'}(X_{i})}{p_{\theta}(X_{i})}}}\quad \text{with}\quad \psi(u)=\frac{u-1}{u+1}\quad \text{for $u\in[0,+\infty]$}
\]
and the conventions $\psi(+\infty)=1$, $0/0=1$ and $a/0=+\infty$ for every $a>0$. Then, given two parameters $\lambda>0$ and $\beta\in (0,1)$, we define for $\theta\in \Theta$,   
\begin{align}
\gT(\bsX,\theta)
&=\int_{\Theta}\gT(\bsX,\theta,\theta')\frac{\exp\cro{\lambda\gT(\bsX,\theta,\theta')}}{\int_{\Theta}\exp\cro{ \lambda \gT(\bsX,\theta,\theta')}d\gpi(\theta')}d\gpi(\theta').\label{def-TX}
\end{align}
Our $T$-Bayes posterior distribution $\gpi_{\!\bsX}$ admits a density with respect to $\gpi$ which is given by
\begin{equation}\label{def-piX}
\frac{d\gpi_{\!\bsX}}{d\gpi}(\theta)=\frac{\exp\cro{-\lambda(1-\beta)\gT(\bsX,\theta)}}{\int_{\Theta}\exp\cro{-\lambda(1-\beta)\gT(\bsX,\theta)}d\gpi(\theta)}\quad \text{for $\theta\in \Theta$}.
\end{equation}
A $T$-Bayes estimator $\widehat \theta$ is a random variable with distribution $\gpi_{\!\bsX}$ conditionally on the data $\bsX$.
Throughout this paper, we choose the values $\lambda=0.11$ and $\beta=0.09$ proposed in Baraud and Romon~\citeyearpar{Baraud:2026aa} although other values are possible provided that they satisfy the constraints given in their Proposition~2. We then define $\gamma$ as the numerical constant $\gamma=7.9\times 10^{-4}$. 

\subsection{The $\gpi$-complexity and the critical radius at a parameter $\theta\in \Theta$}
In order to describe the performance of our estimator, we introduce some fundamental quantities as defined in Baraud and Romon~\citeyearpar{Baraud:2026aa} (see also the references therein). They were originally introduced for a general loss function $\gell$ on the parameter space $\Theta$ but, in the present paper, we use them for the specific loss defined by $\gell(\theta,\theta')=nh^{2}(P_{\theta},P_{\theta'})$ for $\theta,\theta'\in \Theta$. Given $r\ge 0$, we denote by $\sB(\theta,r)$ the closed ball, relative to the loss $\gell$, centred at $\theta$ with radius $r$, that is,  $\sB(\theta,r)=\{\theta' \in \Theta,\; \gell\pa{\theta,\theta'}\le r\}$. 

\begin{df}
Given a prior $\gpi$ on $\Theta$, we define the $\gpi$-complexity at $\theta\in \Theta$ as the smallest nonnegative number $r(\gpi,\theta)$ that satisfies
\begin{align}
0<\gpi\pa{\sB(\theta,2r)}\le \exp\pa{\gamma r}\gpi\pa{\sB(\theta,r)}\quad\text{for every }r\ge \ray(\gpi,\theta).
\label{def-piC}
\end{align}
The critical radius $\overline r(\gpi,\theta)$ at $\theta$ is the smallest nonnegative number that satisfies 
\begin{equation}\label{def-CR}
\gpi\pa{\sB(\theta,r)}\ge \exp\pa{-\gamma r}\quad\text{for every }r\ge \overline \ray(\gpi,\theta).
\end{equation}
\end{df}
We refer to Baraud~\citeyearpar{baraud2021robust} and Baraud and Romon~\citeyearpar{Baraud:2026aa}  for more details about these quantities and in particular on how they can be linked to some more classical notions of complexity such as the dimension, the entropy etc. 

Since our loss $\gell$ is here bounded by $n$, $\gpi\pa{\sB(\theta,r)}=1$ as soon as $r\ge n$. Consequently, $\ray(\gpi,\theta)\le n$ for every prior $\gpi$ on $\Theta$. 

We also use the following properties whose proofs are deferred to Section~\ref{pf-prop-Prop}. 
\begin{prop}\label{prop-Prop}
Let $\Theta$ be a measurable space. The following properties hold.  
\begin{listi}
\item\label{Prop-i} For every prior $\gpi$ on $\Theta$, $\ray(\gpi,\theta)\le \overline r(\gpi,\theta)$ for all $\theta\in \Theta$.
\item\label{Prop-ii} Let $(\gpi_{m})_{m\in\cM}$ be an at most countable family of probability distributions on $\Theta$ and $(L_{m})_{m\in\cM}$ nonnegative numbers such that $\sum_{m\in\cM}e^{-L_{m}}=1$. Then, $\gpi=\sum_{m\in\cM}e^{-L_{m}}\gpi_{m}$ satisfies
\begin{equation}\label{eq-pim}
\overline\ray\pa{\gpi,\theta}\le \inf_{m\in\cM}\cro{\overline \ray\pa{\gpi_{m},\theta}+\gamma^{-1}L_{m}}\quad \text{for all $\theta\in \Theta$}.
\end{equation}
\item\label{Prop-iii} For each $i\in\{1,\ldots,k\}$, let $\Theta_{i}$ be a parameter space equipped with a prior $\gpi_{i}$ and a nonnegative loss function $\gell_{i}$. For every $i\in\{1,\ldots,k\}$ and $\theta_{i}\in \Theta_{i}$, let $R_{i}(\gpi_{i},\theta_{i})$ be a nonnegative number such that 
\[
\gpi_{i}\pa{\ac{\theta_{i}'\in \Theta_{i}, \gell_{i}(\theta_{i},\theta_{i}')\le r}}\ge e^{-\gamma r}\quad \text{for all $r\ge R_{i}(\gpi_{i},\theta_{i})$}.
\]
If $\prod_{i=1}^{k}\Theta_{i}\subset \Theta$ and $\gell$ is a nonnegative loss function on $\Theta$ that satisfies 
\begin{equation}\label{eq-hleHi}
\gell(\theta,\theta')\le \sum_{i=1}^{k}\gell_{i}(\theta_{i},\theta_{i}') \quad \text{for all $\theta,\theta'\in \prod_{i=1}^{k}\Theta_{i}$,}
\end{equation}
then 
\begin{equation}\label{eq-hleHic}
\overline \ray\pa{\bigotimes_{i=1}^{k}\gpi_{i},\theta}\le \sum_{i=1}^{k}R_{i}(\gpi_{i},\theta_{i})\quad \text{for every $\theta\in \prod_{i=1}^{k}\Theta_{i}$.}
\end{equation}
\end{listi}
\end{prop}

\subsection{A general result}
The following result was established in Baraud~\citeyearpar{baraud2021robust}. 
\begin{thm}\label{thmCasDet}
Let $\etc{X}$ be independent random variables with marginal distributions $P_{1}\et,\ldots,P_{n}\et$ respectively and average distribution $P\et=n^{-1}\sum_{i=1}^{n}P_{i}\et$. Given a dominated statistical model $\sP=\{P_{\theta}=p_{\theta}\cdot\gmu, \theta\in \Theta\}$ for $P\et$ and a prior $\gpi$ on $\Theta$, a $T$-Bayes estimator $\widehat \theta$ with distribution $\gpi_{\bsX}$ defined by \eref{def-piX} satisfies the following properties. There exists a universal constant $C\ge 1$ such that for every $\xi>0$, with probability at least $1-2e^{-\xi}$, 
\begin{equation}\label{eq-thmCasDet02}
h^{2}\pa{P\et,P_{\widehat \theta}}\le C\inf_{\theta\in \Theta}\cro{h^{2}(P\et,P_{\theta})+\frac{\ray(\gpi,\theta)+1+\xi}{n}}.
\end{equation}
In particular, 
\begin{equation}\label{eq-thmCasDet03}
\E\cro{h^{2}\pa{P\et,P_{\widehat \theta}}}\le C\inf_{\theta\in \Theta}\cro{h^{2}(P\et,P_{\theta})+\frac{\ray(\gpi,\theta)+1}{n}}.
\end{equation}
\end{thm}
As explained in the introduction, we now consider the situation where $\bsX$ is an $n$-sample of determinantal processes with distribution $P\et$.

\section{Estimating the distribution of a determinantal projection process}\label{sect-PDP}
In this section, we assume that $P\et$ is the distribution of a $k$-Determinantal Projection Process. We denote by $\sD_{k}$ the set of $k$-PDPs on $(\cX,\cA,\mu)$, that is, the set of probability measures with densities of the form \eref{def-dens} with respect to $\gmu$.

\subsection{Parameterisation of $\sD_{k}$}
Following our introduction, it is convenient to parameterise $\sD_{k}$ by $k$-dimensional subspaces $U$ of $\H$.  Denoting by $\cG_{k}(\H)$ the set of all subspaces of $\H$ of dimension $k\ge 1$, we may therefore write $\sD_{k}=\{P_{U},\; U\in \cG_{k}(\H)\}$. This parameterisation is, however, not one to one. The lack of identifiability of the parameter $U$ can already be seen when $k=1$. If $U$ is  the linear span of a unit vector $u\in\H$ and $V$ that of $|u|$, the densities $p_{U}=p_{V}=|u|^{2}$ coincide while $U$ and $V$ are different in general. This example also shows that the kernel $K(x,y)=u(x) u(y)$ which is associated with the distribution $P_{U}$ is not identifiable either.

Despite this lack of identifiability, the parameterisation of $\sD_{k}$ by $\cG_{k}(\H)$ offers several advantages, one being that it is possible to equip $\cG_{k}(\H)$ with a distance $\frd$ for which $\frd(U,V)$ with $U,V\in\cG_{k}(\H)$ nicely controls the Hellinger distance between $P_{U}$ and $P_{V}$. This distance and its geometrical interpretation is given in the following section. Given  $\frd$, we may equip $\cG_{k}(\H)$ with its Borel $\sigma$-algebra and turn it into a measurable space, which we shall do.   

\subsection{The Hellinger distance on $\sD_{k}$}
Let us first define the principal angles and vectors between two finite-dimensional linear spaces sharing the same dimension. 
\begin{df}
Let $k$ be a positive integer and $U$ and $V$ two $k$-dimensional linear subspaces of a Hilbert space $\H$ equipped with its inner product $\scal{\cdot}{\cdot}$. We say that $0\le \fra_{1}\le \ldots\le \fra_{k}\le \pi/2$ are {\em the principal angles} between the linear spaces $U$ and $V$ if there exist an orthonormal basis $\{u_{1}\et,\ldots,u_{k}\et\}$ of $U$ and an orthonormal basis $\{v_{1}\et,\ldots,v_{k}\et\}$ of $V$ such that for $i,j\in\{1,\ldots,k\}$
\[
\scal{u_{i}\et}{v_{i}\et}=\cos \fra_{i}\quad \text{while for $i\ne j$}\quad \scal{u_{i}\et}{v_{j}\et}=0.
\]
The vectors $\{u_{1}\et,\ldots,u_{k}\et\}$, $\{v_{1}\et,\ldots,v_{k}\et\}$ are called {\em principal vectors} of $U$ and $V$.
\end{df}
Given $U,V\in \cG_{k}(\H)$, we denote by $\fra_{1}(U,V)\le \ldots\le \fra_{k}(U,V)$ the principal angles between $U$ and $V$ sorted in increasing order. Given an orthonormal basis $\gu=\{u_{1},\ldots,u_{k}\}$ of $U$ and an orthonormal basis $\gv=\{v_{1},\ldots,v_{k}\}$ of $V$, the cosines of the principal angles between $U$ and $V$ are the singular values of the matrix $G(\gu,\gv)=\pa{\scal{u_{j}}{v_{j'}}}_{1\le j,j'\le k}$. The principal angles can therefore be obtained by finding the singular values of $G(\gu,\gv)$, the result being independent of the choice of the orthonormal bases $\gu$ and $\gv$ for $U$ and $V$.  In particular, the product
\begin{equation}\label{eq-def-d0}
\prod_{i=1}^{k}\cos\cro{\fra_{i}(U,V)}=\ab{\det G(\gu,\gv)}
\end{equation}
is independent of $\gu$ and $\gv$. It is equal to 1 if and only if $U=V$. The following result is proven in Section~\ref{pf-prop-DPP}.
\begin{prop}\label{prop-DPP}
Let $k\ge 1$. The function $\frd$ on $\cG_{k}(\H)\times\cG_{k}(\H)$ defined by 
\begin{equation}\label{eq-def-d}
\frd(U,V)=\sqrt{1-\prod_{i=1}^{k}\cos\cro{\fra_{i}(U,V)}},
\end{equation}
is a distance on $\cG_{k}(\H)$. Besides, if  $U,V\in\cG_{k}(\H)$, 
\begin{equation}\label{eq-DPP-aff}
h^{2}(P_{U},P_{V})\le \frd^{2}(U,V)\le \frac{1}{2}\sum_{i=1}^{k}\norm{u_{i}-v_{i}}^{2}
\end{equation}
where $\{u_{1},\ldots,u_{k}\}$ and $\{v_{1},\ldots,v_{k}\}$ are arbitrary orthonormal bases of $U$ and $V$ respectively. 
\end{prop}
The constant $1/2$ in the right-hand side of \eref{eq-DPP-aff} improves that obtained in Baraud~\citeyearpar{Baraud:2012qy} [Proposition 3.4]. If $u$ and $v$ are two nonnegative functions in the unit sphere of $\H$ and $U,V$ their respective linear spans,  $P_{U}=u^{2}\cdot\mu$, $P_{V}=v^{2}\cdot\mu$,   
\begin{equation}\label{eq-h-L2}
h^{2}(P_{U},P_{V})=\frac{1}{2}\int_{\cX}\pa{u-v}^{2}d\mu=\frac{1}{2}\norm{u-v}^{2}
\end{equation}
which shows that the constant $1/2$ in \eref{eq-DPP-aff} cannot be improved.

\subsection{The Haar measure $\frm_{k}^{E}$ on the Stiefel manifold}\label{sect-StM}
Given $k\ge 1$ and a Euclidean space $E$ with finite dimension $D\ge k$, the Stiefel manifold $\U_{k}(E)$ is the set of $k$-tuples $\gu=(u_{1},\ldots,u_{k})\in E^{k}$ for which $u_{1},\ldots,u_{k}$ is an orthonormal system of $E$. Given an orthonormal basis $e_{1},\ldots,e_{D}$ of $E$, we may write these elements in the form $\gu_{H}=(\sum_{i=1}^{D}H_{i,1}e_{i},\ldots, \sum_{i=1}^{D}H_{i,k} e_{i})$ where the columns of the $D\times k$ matrix $H=(H_{i,j})$ 
 are orthonormal in $\R^{D}$, that is, satisfy  $(H^{t})H=I_{k}$ where $I_{k}$ denotes the  $k\times k$ identity matrix. There exists a unique probability measure $\frm_{k}^{E}$ on  $\U_{k}(E)$ that remains invariant under the mapping   $\gu\mapsto O\gu=(Ou_{1},\ldots,Ou_{k})$ for every orthogonal transformation $O$ of $E$. It is called the {\em Haar measure on the Stiefel manifold $\U_{k}(E)$.} We refer the reader to Muirhead~\citeyearpar{MR652932} for more details. When $k=1$,  $\frm_{1}^{E}$ is simply the uniform distribution on the unit sphere of $E$. 
  
 We can simulate a random variable $\gu$ with distribution $\frm_{k}^{E}$ as follows. First generate a $D\times k$ matrix $Z$ whose entries are i.i.d.\ standard Gaussian random variables. The matrix $Z$ a.s.\ admits a unique decomposition $Z=HT$ where $H$ is a $D\times k$ matrix whose columns are orthonormal and $T$ a $k\times k$ upper-triangular matrix with positive elements on the diagonal. This decomposition corresponds to a Gram--Schmidt orthogonalisation procedure applied to the columns of $Z$. Since the distribution of $Z$ is invariant under orthogonal transformation, $\gu_{H}=(\sum_{i=1}^{D}H_{i,1}e_{i},\ldots, \sum_{i=1}^{D}H_{i,k} e_{i})$ has the required distribution $\frm_{k}^{E}$. We use this probability measure in order to define our prior distribution. 

\subsection{A prior on $\cG_{k}(\H)$}\label{sect-prior01}
In this section, our aim is to propose a prior $\gpi$ on the parameter space $\Theta=\cG_{k}(\H)$ which yields a posterior distribution (or a randomised estimator) with good properties. 

Let $E$ be a linear subspace of $\H$ of dimension $D>k$. If we believe that the true distribution of the data is of the form $P_{U}$ for some $k$-dimensional subspace $U$ of $E$, we should define our prior on the Grassmannian $\cG_{k}(E)$. To do so, we proceed as follows. We may generate a linear subspace of $E$ of dimension $k$ by considering the linear span ${\rm Span}\ac{v_{1},\ldots,v_{k}}$ of $v_{1},\ldots,v_{k}$ when $\gv=(v_{1},\ldots,v_{k})$ has the Haar distribution $\frm_{k}^{E}$ on the Stiefel manifold $\U_{k}(E)$. The image $\gpi$  of $\frm_{k}^{E}$ by the mapping $\gv=(v_{1},\ldots,v_{k})\mapsto {\rm Span}\ac{v_{1},\ldots,v_{k}}$ yields a distribution on the Grassmannian $\cG_{k}(E)$. Since $\frm_{k}^{E}$ remains invariant under orthogonal transformations of $E$, $\gpi$ is also invariant under the mapping $V\in\cG_{k}(E)\mapsto OV=\{Ov,\; v\in V\}$ for every orthogonal transformation $O$ of $E$. This probability is the {\em Haar measure} on the Grassmannian $\cG_{k}(E)$ and it satisfies the following property. The result is proven in Section~\ref{pf-prop-rUni}.
\begin{prop}\label{prop-rUni}
Let $k\ge 1$ and $E$ be a linear subspace of $\H$ with finite dimension $D\ge k+1$. Let us set $a=\min\{k,D-k\}$ and $b=\max\{k,D-k\}$. The Haar measure $\gpi$ on the Grassmannian $\cG_{k}(E)$ satisfies for every $U\in\cG_{k}(E)$ and $t\in [0,1]$, 
\begin{align}
\gpi\pa{\ac{V\in\cG_{k}(E), \frd(U,V)\le t}}&=\P\cro{\prod_{i=1}^{a}B_{i}\ge \pa{1-t^{2}}^{2}}\label{eq-rUni00}
\end{align}
where the $B_{i}$ are independent random variables with Beta distributions $\cB(i/2,b/2)$ respectively for $i\in\{1,\ldots,a\}$. 
Moreover, 
\begin{equation}\label{eq-rUni01}
\gpi\pa{\ac{V\in\cG_{k}(E), \frd(U,V)\le t}}\ge \pa{\frac{t}{2\sqrt{k\wedge (D-k)}}}^{k(D-k)}\quad \text{for every $t\in [0,1]$.}
\end{equation}
\end{prop}
When $\dim E=D=k$, $\gpi$ is the Dirac mass $\delta_{E}$ at $E$ and  \eref{eq-rUni01} still holds with the convention $\infty^{0}=1$. We use this convention hereafter. 

We immediately derive the following corollary whose proof is deferred to Section~\ref{pf-cor-cont01}. 
\begin{cor}\label{cor-cont01}
Let $k\ge 1$ and $E$ be a linear subspace of $\H$ with finite dimension $D\ge k$. If we equip the statistical model $\{P_{U},\; U\in\cG_{k}(\H)\}$ with the Haar probability measure $\gpi$ on the Grassmannian $\cG_{k}(E)$, we obtain that
\begin{equation}\label{eq-rayGk}
\overline \ray(\gpi,U)\le \frac{k(D-k)}{2\gamma}\log_{+}\pa{\frac{8\gamma n }{k\vee (D-k)}}\quad \text{for every $U\in\cG_{k}(E)$}.
\end{equation}
\end{cor}
Dealing with a single linear subspace $E$ of $\H$ might be too restrictive and it may be wiser to consider an at most countable collection $\{E_{m},\; m\in\cM\}$ of candidate ones. Denoting by $\gpi_{m}$ the Haar probability measure on the Grassmannian $\cG_{k}(E_{m})$ and nonnegative numbers $\{L_{m},\; m\in\cM\}$ that satisfies $\sum_{m\in\cM}e^{-L_{m}}=1$, we may consider the prior $\gpi$ on $\cG_{k}(\H)$ defined by 
\begin{equation}\label{def-pimG}
\gpi=\sum_{m\in\cM}e^{-L_{m}}\gpi_{m}.
\end{equation}
It then follows from \eref{eq-pim} that
\begin{equation}\label{eq-bornerayU}
\overline \ray(\gpi,U)\le \inf_{m\in\cM}\cro{\overline \ray(\gpi_{m},U)+\gamma^{-1}L_{m}}\quad \text{for every $U\in\cG_{k}(\H)$.}
\end{equation}
The following result is a straightforward consequence of Theorem~\ref{thmCasDet} applied to the statistical model $\{P_{U},\; U\in \cG_{k}(\H)\}$. 
\begin{cor}\label{cor-casPDPex}
Let $\{E_{m},\; m\in\cM\}$ be an at most countable family of finite-dimensional linear subspaces of $\H$ with respective dimensions $D_{m}\ge k\ge 1$. Equip the parameter space $\cG_{k}(\H)$ with the prior $\gpi$ defined by \eref{def-pimG}. The $T$-Bayes estimator $\widehat U$ with distribution $\gpi_{\bsX}$ on $\cG_{k}(\H)$ satisfies the following property. If $\etc{X}$ are i.i.d.\ with distribution $P_{U}$ for some $U\in \cG_{k}(E_{m})$ and $m\in\cM$, then for every $\xi>0$, with probability at least $1-2e^{-\xi}$, 
\begin{equation}\label{eq-casPDPex}
h^{2}\pa{P_{U},P_{\smash{\widehat U}}}\le C\cro{\frac{k(D_{m}-k)}{n}\log_{+}\pa{\frac{n}{k\vee(D_{m}-k)}} +\frac{L_{m}+1 +\xi}{n}}.
\end{equation}
\end{cor}

If the collection $\{E_{m},\; m\in\cM\}$ restricts to a singleton $E$ of dimension $D\ge k$, we may take $L_{m}=0$,  and by integrating \eref{eq-casPDPex} with respect to $\xi$ we obtain that for $U\in\cG_{k}(E)$ the risk of $\widehat U$ satisfies 
\begin{equation}\label{eq-casPDPex1}
\E\cro{h^{2}\pa{P_{U},P_{\smash{\widehat U}}}}\le  C\cro{\frac{k(D-k)}{n}\log_{+}\pa{\frac{n}{k\vee(D-k)}} +\frac{1}{n}}.
\end{equation}
When $D>k$, the rate of convergence of $P_{\widehat U}$ to $P_{U}$ for the Hellinger loss is therefore at most $\sqrt{(k(D-k)\log n)/n}$. 

If one considers now a collection $\{E_{m},\; m\in\cM\}$ of such linear spaces and $U$ belongs to one of them, say $\cG_{k}(E_{m})$, the risk bound is inflated only by the additional term $L_{m}/n$. 


\subsection{Approximating $U$}
Assuming that the data are i.i.d.\ with a distribution of the form $P_{U}$ for some  $U$ in $\cG_{k}(E)$ (or in $\bigcup_{m\in\cM}\cG_{k}(E_{m})$), as in the previous section, is too strong in general. A more reasonable approach is to assume that $P_{U}$ is well approximated (in the Hellinger distance) by an element of the form $P_{V}$ with $V$ in $\cG_{k}(E)$. A look at \eref{eq-thmCasDet02} shows that compared to the ideal situation where $P\et=P_{U}$ with $U\in \cG_{k}(E)$,  assuming that $P_{U}$ is close to a distribution $P_{V}$ with $V\in\cG_{k}(E)$ can only inflate the right-hand side of  \eref{eq-casPDPex}
by an extra term of order $h^{2}(P_{U},P_{V})$. This raises the following question. How can we control the quantity $\sA(U,E)=\inf_{V\in\cG_{k}(E)}h\pa{P_{U},P_{V}}$? If the linear space $U$ is generated by the orthonormal basis $\{u_{1},\ldots,u_{k}\}$, is it possible to link the quantity $\sA(U,E)$ to the way the linear space $E$ approximates each element of the orthonormal basis $u_{1},\ldots,u_{k}$? The answer to this question is given by the following result proven in Section~\ref{pf-prop-approx01}
\begin{prop}\label{prop-approx01}
Let $k\ge 1$ and $U$ be an element of $\cG_{k}(\H)$ generated by the orthonormal basis $\{u_{1},\ldots,u_{k}\}$. If $E$ is a finite-dimensional linear subspace of $\H$ of dimension at least $k$,  
\begin{equation}\label{eq-DPP-aff01}
 \inf_{V\in\cG_{k}(E)}h^{2}(P_{U},P_{V})\le \inf_{V\in\cG_{k}(E)}\frd^{2}(U,V)\le \sum_{i=1}^{k}\inf_{v\in E}\norm{u_{i}-v}^{2}.
\end{equation}
\end{prop}
Note that the second inequality in \eref{eq-DPP-aff01} becomes an equality when $k=1$.

Another natural question is how to choose $E$. In order to provide an example, let us assume for a moment that $\H$ consists of functions on an interval $\cX$ of $\R$ equipped with the Lebesgue measure $\mu$. The kernel associated with $P_{U}$ is then $K(x,y)=\sum_{i=1}^{k}u_{i}(x)u_{i}(y)$ where the $u_{i}$ satisfy
\[
u_{i}(x)=\int_{\cX}K(x,y)u_{i}(y)dy\quad \text{for $x\in\cX$ and $i\in\{1,\ldots,k\}$.}
\]
If $K$ is smooth with respect to $x\in \cX$, so are the functions $u_{1},\ldots,u_{k}$ on $\cX$. It is well known that smooth functions can be well approximated by finite-dimensional linear spaces (generated by splines, trigonometric polynomials, wavelets, etc.). In view of \eref{eq-DPP-aff01} it is therefore natural to consider  linear subspaces $E$ of this form.

Using Proposition~\ref{prop-approx01}, we can improve the result of Corollary~\ref{cor-casPDPex} by relaxing the assumption that $P_{U}$ belongs to $\bigcup_{m\in\cM}\cG_{k}(E_{m})$. We derive from Theorem~\ref{thmCasDet}, more precisely from \eref{eq-thmCasDet02} the following result.
\begin{cor}\label{cor-casPDPex01}
Let $\{E_{m},\; m\in\cM\}$ be an at most countable family of finite-dimensional linear spaces of $\H$ with respective dimensions $D_{m}\ge k\ge 1$ and $\gpi$ be the prior defined by \eref{def-pimG}. The $T$-Bayes estimator $\widehat U$ with distribution $\gpi_{\bsX}$ on $\cG_{k}(\H)$ satisfies the following property. If $U$ is  generated by the orthonormal basis $\{u_{1},\ldots,u_{k}\}$ and $\etc{X}$ are i.i.d.\ with distribution $P_{U}$, then for every $\xi>0$, with probability at least $1-2e^{-\xi}$, 
\begin{equation}\label{eq-casPDPex2}
h^{2}\pa{P_{\smash{U}},P_{\smash{\widehat U}}}\le C\inf_{m\in\cM}\cro{\cO_{m}(U) +\frac{L_{m}+1 +\xi}{n}}
\end{equation}
where 
\[
\cO_{m}(U)=\sum_{i=1}^{k}\cro{\inf_{v\in E_{m}}\norm{u_{i}-v}^{2}+ \frac{D_{m}\log_{+}\pa{n/D_{m}}}{n}}.
\]
\end{cor}
By integrating \eref{eq-casPDPex2} with respect to $\xi$, we get the risk bound 
\begin{equation}\label{eq-RB01}
\E\cro{h^{2}\pa{P_{\smash{U}},P_{\smash{\widehat U}}}}\le C\inf_{m\in\cM}\cro{\sum_{i=1}^{k}\pa{\inf_{v\in E_{m}}\norm{u_{i}-v}^{2}+\frac{D_{m}}{n}\log_{+} \pa{\frac{n}{D_{m}}}}+\frac{L_{m}}{n}}.
\end{equation}
 If $k=1$ and $U$ is the linear span of $u\ge 0$ in the unit sphere of $\H$, the problem of estimating $P_{U}$ is equivalent to that of estimating the density $p=u^{2}$ from an $n$-sample. Up to the extra logarithmic factor, our risk bound \eref{eq-RB01} is the one we would expect for estimating the density $p$ from the collection of linear spaces $\{E_{m},\; m\in\cM\}$ under the Hellinger loss (since $h$ and the $\L_{2}$-distance are linked by \eref{eq-h-L2}). For $k>1$, our risk bound is not larger than the sum of these individual risk bounds for estimating each of the densities $p_{i}=u_{i}^{2}$ separately.

\section{General determinantal processes}\label{sect-DPP}
In the remainder of the paper, it is convenient to associate DPPs with kernels $K$ of the form 
\begin{equation}\label{eq-Kgene}
K(x,y)=\sum_{j\in\frJ}l_{j}u_{j}(x)u_{j}(y)\quad \text{for $(x,y)\in\cX^{2}$}
\end{equation}
where $\frJ$ is a finite or countable set, that is, to allow $\frJ$ to be different from $\N\et$ as in \eref{def-K}. We recall that the elements of $\gu=(u_{j})_{j\in\frJ}$ are orthonormal in $\H$ and $\gl=(l_{j})_{j\in \frJ}$ is a sequence of numbers in  $[0,1]$ that satisfies $0\le \sum_{j\in \frJ}l_{j}<+\infty$. Of particular interest is the case where $\cX$ is finite, in which case we may assume that $\cX=\{1,\ldots,N\}$, $\H=\R^{N}$ and $\frJ=\cX$. 

We denote by $\sD$ the class of all DPP distributions on $\cX$ associated with kernels of the form~\eref{eq-Kgene}. In order to parameterise $\sD$, we introduce the set $\U$ of the sequences $\gu=(u_{j})_{j\in\frJ}$ whose elements are orthonormal in  $\H$ and $\gS$ that of the sequences $\gs=(s_{j})_{j\in\frJ}$ that belong to $[0,\pi/2]^{\frJ}$ and satisfy $\sum_{j\in\frJ} s_{j}^{2}<+\infty$. The class of all finite subsets $J$  of $\frJ$ is denoted by $\cJ$. Given $\gu\in \U$ and $J\in\cJ\setminus\{\vide\}$, it is convenient to use the uppercase letter $U$ together with the subscript $J$ to denote the linear span  $U_{J}$ of the $u_{j}$ for $j\in J$. In the same manner, $V_{J}$ denotes the linear span of the $v_{j}$ for $j\in J$ when $\gv=(v_{j})_{j\in\frJ}$ belongs to $\U$. Through the change of variables $l_{j}=\sin^{2}s_{j}$ for $j\in \frJ$, every sequence $\gl=(l_{j})_{j\in\frJ}$ of numbers in $[0,1]$ satisfying $\sum_{j\in\frJ}l_{j}<+\infty$ can be associated in a unique way with a sequence $\gs\in\gS$. A distribution of $\sD$ can then be written as a (possibly  infinite)  mixture of PDP as follows:
\begin{equation}\label{def-Psu}
P_{\gs,\gu}=\sum_{J\in\cJ}p_{J}(\gs)P_{U_{J}}\quad \text{with}\quad p_{J}(\gs)=\prod_{j\in J}\sin^{2}\pa{s_{j}}\prod_{j\not\in J}\cos^{2}(s_{j})\quad \text{for}\quad \text{$(\gs,\gu)\in \gS\times \U$}
\end{equation}
with the convention $P_{U_{\vide}}=\delta_{\vide}$.  Under our condition $\sum_{j\in\frJ} s_{j}^{2}<+\infty$, $\prod_{j\in J}\sin^{2}\pa{s_{j}}=0$ for every infinite subset $J$ of $\frJ$ and the equality 
\[
1=\prod_{j\in \frJ}\pa{\cos^{2}\pa{s_{j}}+\sin^{2}\pa{s_{j}}}=\sum_{J\in \cJ}p_{J}(\gs)
\]
is therefore satisfied, showing that $J\mapsto p_{J}(\gs)$ defines a genuine probability measure on $\cJ$ for every parameter $\gs\in\gS$. When $I=\{j\in\frJ,\; s_{j}\ne 0\}$ is nonempty and finite, the sum in \eref{def-Psu} may be restricted to the subsets $J$ of $I$. The distribution $P_{\gs,\gu}$ then only depends on $\gs$ and $\gu_{I}=(u_{j})_{j\in I}$. We then write this distribution as $P_{\gs,\gu_{I}}$.

\subsection{The Hellinger distance on $\sD$}
We have seen that the elements of $\sD$ can be parameterised by a pair $(\gs,\gu)$ where the elements of $\gs$ belong to $[0,\pi/2]$ and are square-summable while $\gu\in \U$ is a sequence of orthonormal vectors of $\H$. For $J\in\cJ\setminus\{\vide\}$, let us define the set 
\[
\gS_{J}=\ac{\gs\in \gS,\; s_{j}=0\text{ for all $j\not\in J$ }},
\]
which we may identify with $[0,\pi/2]^{J}$. For $\gs\in \gS_{J}$, the distribution $P_{\gs,\gu}$ may then be written as $P_{\gs,\gu_{J}}$ where $\gu_{J}=(u_{j})_{j\in J}$ is a finite orthonormal system of $\H$. If these vectors belong to a linear space $E$ of dimension $D\ge |J|$, $\gu_{J}$ is an element of the Stiefel manifold $\U_{J}(E)$. A nice feature of this parameterisation of $\sD$ is the following result. Its proof is deferred to Section~\ref{pf-prop-h-DPP}.
\begin{prop}\label{prop-h-DPP}
Let $(\gs,\gu)\in \gS\times \U$ and $J\in \cJ\setminus\{\vide\}$. For every $(\gt,\gv)\in \gS\times \U$, 
\begin{align}\label{eq-prop-h-DPP0}
h^{2}\pa{P_{\gs,\gu},P_{\gt,\gu}}\le  \frac{1}{2}\sum_{j\in\frJ}\ab{s_{j}-t_{j}}^{2},\quad h^{2}\pa{P_{\gt,\gu},P_{\gt,\gv}}&\le \frac{1}{2}\sum_{j\in \frJ}\sin^{2}(t_{j})\norm{u_{j}-v_{j}}^{2}.
\end{align}
In particular, if $\overline \gs=(s_{j}\1_{j\in J})_{j\in\frJ} \in\gS_{J}$,
\begin{equation}\label{eq-prop-h-DPP01b}
h^{2}\pa{P_{\gs,\gu},P_{\overline\gs,\gu}}\le \frac{1}{2}\sum_{j\not\in J}s_{j}^{2}.
\end{equation}
Furthermore, if $E$ is a linear subspace of $\H$ of dimension $D\ge |J|$, there exists $\overline \gu=(\overline u_{j})_{j\in J}\in\U_{J}(E)$ such that
\begin{align}\label{eq-prop-h-DPP02b}
h^{2}\pa{P_{\gt,\gu},P_{\gt,\overline \gu}}&\le \sum_{j\in J}\inf_{v\in E}\norm{u_{j}-v}^{2}\quad \text{for any $\gt\in\gS_{J}$.}
\end{align}
\end{prop}
Since $h^{2}(P_{\gs,\gu},P_{\gt,\gv})\le 2h^{2}(P_{\gs,\gu},P_{\gt,\gu})+2h^{2}(P_{\gt,\gu},P_{\gt,\gv})$, 
we immediately infer from \eref{eq-prop-h-DPP01b} and \eref{eq-prop-h-DPP02b} that for every $J\in\cJ\setminus\{\vide\}$ there exists $(\overline \gs,\overline \gu)\in\gS_{J}\times \U_{J}(E)$ such that 
\begin{align}\label{eq-prop-h-DPP03b}
h^{2}(P_{\gs,\gu},P_{\overline \gs,\overline\gu})\le 2\cro{\frac{1}{2}\sum_{j\not \in J}s_{j}^{2}+\sum_{j\in J}\inf_{v\in E}\norm{u_{j}-v}^{2}}.
\end{align}

\subsection{Priors on $\sD$ and statistical results}
We have seen that given $J\in\cJ\setminus\{\vide\}$, a DPP distribution in $\sD$ can be approximated by elements of the form $P_{\gt,\gv}$ with $(\gt,\gv)\in \gS_{J}\times \U_{J}(E)$ where $\U_{J}(E)$ is the Stiefel manifold associated with a linear subspace $E$ of $\H$ of dimension $D\ge |J|$. This suggests considering priors supported on subsets of the parameter space $\gS\times \U$ which are of the form $\gS_{J}\times\U_{J}(E)$. 

A consequence of Proposition~\ref{prop-h-DPP} is that the restriction of our loss $\gell$ to $\gS_{J}\times \U_{J}(E)$ satisfies the following property. For $(\gs,\gu)$ and $(\gt,\gv)$ in $\gS_{J}\times \U_{J}(E)$, 
\begin{align}\label{eq-ellmaj}
\gell((\gs,\gu),(\gt,\gv))=nh^{2}(P_{(\gs,\gu)},P_{(\gt,\gv)})\le \gell_{1}(\gs,\gt)+\gell_{2}(\gu,\gv)
\end{align}
with $\gell_{1}(\gs,\gt)=n\sum_{j\in J}\ab{s_{j}-t_{j}}^{2}$ and $\gell_{2}(\gu,\gv)=n\sum_{j\in J}\norm{u_{j}-v_{j}}^{2}$. Applying Proposition~\ref{prop-Prop} we derive the following result proven in Section~\ref{pf-prop-produit}.
\begin{prop}\label{prop-produit}
Let $J\in\cJ\setminus\{\vide\}$, $E$ be a linear subspace of $\H$ of dimension $D\ge |J|$ and $(\gs,\gu)$ an element of $\gS_{J}\times \U_{J}(E)$. If $\gpi_{1}$ is a prior on $\gS_{J}$ that satisfy

\begin{equation}\label{def-produit1}
\gpi_{1}\pa{\ac{\gt\in \gS_{J}, \sum_{j\in J}\ab{s_{j}-t_{j}}^{2}\le \frac{r}{n}}}\ge e^{-\gamma r}\quad \text{for all $r\ge R_{1}(\gpi_{1},\gs)$}
\end{equation}
and $\gpi_{2}$ is a prior on the Stiefel manifold $\U_{J}(E)$ that satisfies
\begin{equation}\label{def-produit2}
\gpi_{2}\pa{\ac{\gv\in \U_{J}(E), \sum_{j\in J}\norm{u_{j}-v_{j}}^{2}\le \frac{r}{n}}}\ge e^{-\gamma r}\quad \text{for all $r\ge R_{2}(\gpi_{2},\gu)$}
\end{equation}
then the prior $\gpi=\gpi_{1}\otimes\gpi_{2}$ on $\gS_{J}\times\U_{J}(E)\subset \gS\times \U$ satisfies
\begin{equation}\label{eq-prop-produit}
\overline \ray\pa{\gpi,(\gs,\gu)}\le R_{1}(\gpi_{1},\gs)+R_{2}(\gpi_{2},\gu).
\end{equation}
In particular, for  the uniform distribution $\gpi_{1}=\cU(\gS_{J})$ on $\gS_{J}=[0,\pi/2]^{J}$ we may take
\begin{equation}\label{eq-prop-produit1}
R_{1}(\gpi_{1},\gs)=\frac{|J|}{2\gamma }\log_{+}\pa{\frac{\gamma n\pi^{2}}{2}}
\end{equation}
and for the Haar measure $\gpi_{2}=\frm_{J}^{E}$ on $\U_{J}(E)$,
\begin{equation}\label{eq-prop-produit2}
R_{2}(\gpi_{2},\gu)=\frac{|J|D}{2\gamma}\log\cro{4\pi\pa{\frac{4\gamma n}{D}\vee 1}}.
\end{equation}

\end{prop}

Given an $n$-sample $\etc{X}$ with unknown distribution $P_{\gs,\gu}\in\sD$, we may apply our $T$-Bayes estimation procedure described in Section~\ref{sect-TPost} to estimate $P_{\gs,\gu}$ by using a prior of the form $\gpi=\gpi_{1}\otimes \gpi_{2}$ where $\gpi_{1}$ and $\gpi_{2}$ satisfy the assumptions of Proposition~\ref{prop-produit}. Using \eref{eq-prop-h-DPP0}, Proposition~\ref{prop-Prop}-\ref{Prop-i} and \eref{eq-prop-produit} and applying inequality \eref{eq-thmCasDet02} 
to the probability $P\et=P_{\gs,\gu}$ we obtain
\begin{cor}\label{cor-su}
Let $\etc{X}$ be an $n$-sample with distribution $P_{\gs,\gu}\in\sD$, $J$ an element of $\cJ\setminus\{\vide\}$ and $E$ a linear subspace of $\H$ of dimension $D\ge |J|$. If $\gpi$ is a prior on $\Theta=\gS\times\U$ of the form $\gpi=\gpi_{1}\otimes \gpi_{2}$ where $\gpi_{1}$ and $\gpi_{2}$ satisfies the assumptions of Proposition~\ref{prop-produit} then for every $(\gt,\gv)\in \gS_{J}\times \U_{J}(E)$ and $\xi>0$, with probability at least $1-2e^{-\xi}$, a $T$-Bayes estimator $\widehat \theta=(\widehat \gs,\widehat \gu)$ with distribution $\gpi_{\bsX}$ satisfies 
\begin{align}
C^{-1}h^{2}\pa{P_{\gs,\gu},P_{\smash{\widehat \gs,\widehat \gu}}}&\le h^{2}(P_{\gs,\gu},P_{\gt,\gv})+\frac{\ray(\gpi,(\gt,\gv))+1+\xi}{n}\nonumber\\
&\le  h^{2}(P_{\gs,\gu},P_{\gt,\gv}) +\frac{ R_{1}(\gpi_{1},\gt)+R_{2}(\gpi_{2},\gv)+1+\xi}{n}.\label{eq-su0}
\end{align}
In particular, by choosing  $(\gt,\gv)=(\overline \gs,\overline\gu)\in\gS_{J}\times\U_{J}(E)$ that satisfies \eref{eq-prop-h-DPP03b}, we obtain that $(\widehat \gs,\widehat \gu)$ satisfies
\begin{align}\label{eq-su}
C^{-1}h^{2}\pa{P_{\gs,\gu},P_{\smash{\widehat \gs,\widehat \gu}}}&\le \sum_{j\not \in J}s_{j}^{2}+\sum_{j\in J}\inf_{v\in E}\|u_{j}-v\|^{2}+\frac{ R_{1}(\gpi_{1},\overline \gs)+R_{2}(\gpi_{2},\overline \gu)+1+\xi}{n}.
\end{align}
\end{cor}

Applying \eref{eq-su} to $\gpi=\cU(\gS_{J})\otimes \frm_{J}^{E}$ and integrating with respect to $\xi$, we derive from \eref{eq-prop-produit1} and \eref{eq-prop-produit2} the risk bound 
\begin{align}
C^{-1}\E\cro{h^{2}\pa{P_{\gs,\gu},P_{\smash{\widehat \gs,\widehat \gu}}}}\le \sum_{j\not \in J}s_{j}^{2}+\sum_{j\in J}\cro{\inf_{v\in E}\norm{u_{j}-v}^{2}+\frac{D}{n}\log_{+}\pa{\frac{n}{D}}}.\label{dev-1model}
\end{align}
Let us now analyse this  bound. The term $\sum_{j\not \in J}s_{j}^{2}$ is the error that we make for keeping in \eref{eq-Kgene} only the indices $j$ that belong to the finite set $J$. This amounts to approximating the kernel $K$ defined by \eref{eq-Kgene} by 
\[
K_{J}(x,y)=\sum_{j\in J}l_{j}u_{j}(x)u_{j}(y)\quad \text{for $(x,y)\in\cX^{2}$.}
\]
For $j\in J$, the term
\[
\inf_{v\in E}\norm{u_{j}-v}^{2}+\frac{D}{n}\log_{+}\pa{\frac{n}{D}}
\]
is, up to the logarithmic factor, the bound we would get for estimating the density $p=u_{j}^{2}$ alone from the statistical model $\{P_{u}=u^{2}\cdot\mu,\; u\in \S_{E}\}$ where $\S_{E}$ denotes the unit sphere of $E$. The sum over $j$ in $J$ of these bounds therefore corresponds to the estimation error that we get for estimating all the $u_{j}$ for $j\in J$ simultaneously.

Considering a single finite subset $J$ of $\frJ$ and a single linear subspace $E$ of $\H$ might be too restrictive. We may also consider a family $\cJ'$ of possible elements $J\in\cJ\setminus\{\vide\}$ and associate each of them with an at most  countable collection $\{E_{m},\; m\in\cM_{J}\}$ of linear subspaces $E_{m}$ of $\H$ of dimensions $D_{m}\ge |J|$. For this purpose, we introduce a prior $\nu$ on the set $\frM=\{(J,m),\; J\in \cJ', m\in\cM_{J}\}$ as well as a collection of product probabilities $\gpi_{1,J}\otimes \gpi_{2,m}$ on $\gS_{J}\times \U_{J}(E_{m})$ for $(J,m)\in \frM$ in order to equip $\gS\times\U$ with the hierarchical prior
\begin{equation}\label{def-hprior}
\gpi=\sum_{(J,m)\in\frM}\nu(J,m)\pa{\gpi_{1,J}\otimes\gpi_{2,m}}.
\end{equation}
We immediately derive from \eref{eq-thmCasDet02} and Proposition~\ref{prop-Prop}-\ref{Prop-ii} the following result.
\begin{cor}\label{cor-su01}
For every $(J,m)\in\frM$, let $(\overline \gs_{J},\overline \gu_{J,m})$ be an element of $\gS_{J}\times\U_{J}(E_{m})$
that satisfies \eref{eq-prop-h-DPP03b} with $E=E_{m}$. If $\etc{X}$ is an $n$-sample with distribution $P_{\gs,\gu}\in\sD$ and $\gpi$ the prior on $\Theta=\gS\times\U$ defined by \eref{def-hprior}, then for every $\xi>0$, with probability at least $1-2e^{-\xi}$, a $T$-Bayes estimator $\widehat \theta=(\widehat \gs,\widehat \gu)$ with distribution $\gpi_{\bsX}$ satisfies 
\begin{align}
C^{-1}h^{2}\pa{P_{\gs,\gu},P_{\widehat \gs,\widehat \gu}}\le \inf_{J\in\cJ'}\cro{ \sum_{j\not \in J}s_{j}^{2}+\frac{R_{1}(\gpi_{1,J},\overline \gs_{J})}{n}+\cO_{J}(\gu)+\frac{1+\xi}{n}}\label{eq-thm-DG01}
\end{align}
where 
\[
\cO_{J}(\gu)=\inf_{m\in\cM_{J}}\cro{\pa{\sum_{j\in J}\inf_{v\in E_{m}}\norm{u_{j}-v}^{2}}+\frac{R_{2}(\gpi_{2,m},\overline \gu_{J,m})}{n}-\frac{\log \nu(J,m)}{n}}.
\]
\end{cor}
In particular, by integrating with respect to $\xi$ we get
\begin{align}\label{eq-thm-DG01b}
C^{-1}\E\cro{h^{2}\pa{P_{\gs,\gu},P_{\widehat \gs,\widehat \gu}}}\le \inf_{J\in\cJ'}\cro{\sum_{j\not\in J}s_{j}^{2}+\frac{R_{1}(\gpi_{1,J},\overline \gs_{J})}{n}+\cO_{J}(\gu)}.
\end{align}
Of particular interest is the case where $\gpi_{1,J}$ is the uniform distribution on $\gS_{J}$ and $\gpi_{2,m}$ is the Haar measure $\frm_{J}^{E_{m}}$ on $\U_{J}(E_{m})$ for every $(J,m)\in\frM$. It follows from \eref{eq-prop-produit1} and \eref{eq-prop-produit2} that in this case the quantities $R_{1}(\gpi_{1,J},\overline \gs_{J})$ and $R_{2}(\gpi_{2,m},\overline \gu_{J,m})$ can be bounded independently of $\overline \gs_{J}$ and $\overline \gu_{J,m}$ respectively, leading to the risk bound,
\begin{align}
&C^{-1}\E\cro{h^{2}\pa{P_{\gs,\gu},P_{\widehat \gs,\widehat \gu}}}\nonumber\\
&\le  \inf_{(J,m)\in\frM}\cro{\sum_{j\not \in J}s_{j}^{2}+\sum_{j\in J}\cro{\inf_{v\in E_{m}}\norm{u_{j}-v}^{2}+\frac{D_{m}}{n}\log_{+}\pa{\frac{n}{D_{m}}}}-\frac{\log \nu(J,m)}{n}}.\label{dev-mmodel}
\end{align}
Compared to \eref{dev-1model}, \eref{dev-mmodel} is an oracle-type inequality with respect to the collection of parameter spaces $\{\gS_{J}\times \U_{J}(E_{m}),\; (J,m)\in\frM\}$.

Let us now turn to an example.
\begin{exa}
When $\cX$ is finite, say $\cX=\{1,\ldots,N\}$ with $N\ge 1$, $\H$ is isometric to $\R^{N}$ equipped with the Euclidean norm. A symmetric kernel $K$ on $\cX$ can be seen as an $N\times N$ semi-definite matrix that satisfies \eref{eq-Kgene} with $\frJ=\cX$. The associated DPP distribution on $\cX$, denoted by $P_{K}$, takes the form $P_{K}=P_{\gs,\gu}$ with $\gs=(s_{1},\ldots,s_{N})\in \gS_{\cX}=[0,\pi/2]^{N}$ and $\gu=(u_{1},\ldots,u_{N})\in\U_{\cX}(\H)$. The values of $\sin^{2}(s_{j})$ and the vectors $u_{j}$ are the eigenvalues and the eigenvectors of the matrix $K$. To estimate $P_{K}$, we may apply Corollary~\ref{cor-su} with $J=\cX$ and $E=\H$, which is of dimension $D=N$. We immediately deduce from \eref{eq-prop-produit1}, \eref{eq-prop-produit2} and \eref{eq-su} that for the prior $\gpi=\cU(\gS_{\cX})\otimes \frm_{\cX}^{\H}$, the $T$-Bayes estimator $(\widehat \gs,\widehat \gu)\in [0,\pi/2]^{N}\times  \U_{\cX}(\H)$ satisfies  for every $\xi>0$, with probability at least $1-2e^{-\xi}$, 
\begin{align}\label{eq-su-fini}
C^{-1}h^{2}\pa{P_{\gs,\gu},P_{\smash{\widehat \gs,\widehat \gu}}}&\le \frac{N^{2}\log(n+1)+1+\xi}{n}.
\end{align}
Given $(\widehat \gs,\widehat \gu)$, we may also consider the estimator $\widehat K$ of $K$ given by the $N\times N$ matrix $\widehat K=\sum_{j=1}^{N}\sin^{2}(\widehat s_{j})\widehat u_{j}\widehat u_{j}\tr$. In order to relate $\widehat K$ to $K$ we may use the connection between the Hellinger and the total variation distances together with the fact that if $X\sim P_{K}$ 
\[
\P(J\subset X)=\det K_{J}\quad \text{with}\quad K_{J}=(K_{i,j})_{i,j \in J}\quad \text{for every $J\subset \cX$.}
\]
Denoting by  $\norm{P-Q}_{\rm TV}$ the total variation distance between two probabilities $P$ and $Q$, we deduce from \eref{eq-su-fini} that for every $\xi>0$, with probability at least $1-2e^{-\xi}$, 
\begin{align}
\sup_{J\subset \cX,J\ne \vide}\ab{\det K_{J}-\det \widehat K_{J}}&\le \norm{P_{K}-P_{\smash{\widehat K}}}_{{\rm TV}}\le \sqrt{2}h\pa{P_{\gs,\gu},P_{\smash{\widehat \gs,\widehat \gu}}}\nonumber\\
&\le C\cro{N \sqrt{\frac{\log(n+1)}{n}}+\sqrt{\frac{1+\xi}{n}}}\label{eq-su-fini-K}.
\end{align}
In particular, 
\[
\E\cro{\sup_{J\subset \cX,J\ne \vide}\ab{\det K_{J}-\det \widehat K_{J}}}\le C N \sqrt{\frac{\log(n+1)}{n}}
\]

Our estimator $\widehat K$ allows us to estimate all the principal minors of the matrix $K$ simultaneously with an error of order $N\sqrt{(\log n)/n}$. 
\end{exa}

\section{Translation-invariant kernels on the torus $\T^{d}$}\label{sect-tore}
Throughout this section, we consider the case where $\cX$ is the torus $\T^{d}=\R^{d}/\Z^{d}$.  A function $\overline K$ on the torus $\T^{d}$ can also be viewed as a periodic function on $\R^{d}$ that satisfies $\overline K(x+\gk)=\overline K(x)$ for every $x\in\R^{d}$ and $\gk\in\Z^{d}$ and we systematically identify a function $\overline K$ on the torus $\T^d$ with its periodic version. 

\subsection{Some special functions on the torus}
We denote by $\cK$ the set of functions $\overline K$ on $\T^{d}$ which satisfy the following properties 
\begin{listi}
\item\label{H3} $\overline K$ is continuous on $\T^{d}$ and $\overline K(0)>0$;
\item\label{H1} $\overline K$ is even: for every $u\in \R^{d}$, $\overline K(-u)=\overline K(u)$;
\item\label{H4} the Fourier coefficients $\overline K_{\gk}=\int_{\T^{d}}\overline K(u)\exp\cro{-2i\pi\scal{\gk}{u}}du$ for $\gk\in\Z^{d}$ belong to $[0,1]$;
\item\label{H2} the series $\sum_{\gk\in \Z^{d}}\overline K_{\gk}$ converges.
\end{listi}
An element $\overline K\in\cK$ satisfies the following properties. Firstly, the function $\overline K$ being continuous on the compact set $\T^d$, it is bounded and therefore belongs to $\H=\L_{2}(\T^{d},\mu)$. Under Assumption~\ref{H1} $\overline K_{\gk}$ is a real number and $\overline K_{\gk}=\overline K_{-\gk}$ for every $\gk\in \Z^d$. It is thus convenient to introduce the subset $\Z_{0}^{d}$ of $\Z^{d}$ that consists of the elements $\gk\ne \0$ whose first nonzero coordinate is positive. Then $\Z^{d}$ decomposes into the partition $\Z^d=\{\0\}\cup \Z_{0}^d \cup (-\Z_{0}^d)$ and it follows from  Assumptions~\ref{H4} and \ref{H2} that
\begin{equation}\label{eq-coefF}
0\le \sum_{\gk\in \Z^{d}}\overline K_{\gk}^{2}\le  \sum_{\gk\in \Z^{d}}\overline K_{\gk}=\overline K_{\0}+2\sum_{\gk\in \Z_{0}^{d}}\overline K_{\gk}<+\infty.
\end{equation}
 Since the functions $\frc_{\0}=1$, $\frc_{\gk}(x)=\sqrt{2}\cos(2\pi \scal{\gk}{x})$ and $\frs_{\gk}(x)=\sqrt{2}\sin(2\pi \scal{\gk}{x})$ for $\gk\in\Z_{0}^{d}$ form a Hilbert basis of $\H$,  under Assumption~\ref{H2}, the equality 
%
\begin{equation}\label{Kbar-exp}
\overline K(x)=\sum_{\gk\in\Z^d}\overline K_{\gk}\exp\cro{2i\pi\scal{\gk}{x}}=\overline K_{\0}+2\sum_{\gk\in \Z_{0}^{d}}\overline K_{\gk}\cos(2\pi \scal{\gk}{x})
\end{equation}
holds both pointwise (the series converges uniformly on $\T^{d}$) and in $\H$. In particular, it follows from Parseval's equality and \eref{eq-coefF} that 
\[
\|\overline K\|^{2}=\sum_{\gk\in \Z^{d}}\overline K_{\gk}^{2}\le \overline K_{\0}+2\sum_{\gk\in \Z_{0}^{d}}\overline K_{\gk}=\overline K(0).
\]
Finally, under our assumption $\overline K(0)>0$, the coefficients $\overline  K_{\gk}$ are not all 0 and $\overline K\ne 0$ in $\H$.

\subsection{Special DPPs parameterised by $\cK$ and approximation of their distributions}
Given an element $\overline K\in\cK$, we may consider the kernel $K$ defined for $x,y\in \T^{d}$ by
\begin{align}
K(x,y)&=\overline K(x-y)=\overline K_{\0}+2\sum_{\gk\in \Z_{0}^{d}}\overline K_{\gk}\cos(2\pi \scal{\gk}{x-y})\label{eq-KFour}\\
&=\overline K_{\0}+\sum_{\gk\in \Z_{0}^{d}}\overline K_{\gk}\cro{\frc_{\gk}(x)\frc_{\gk}(y)+\frs_{\gk}(x)\frs_{\gk}(y)}.\nonumber
\end{align}
The values $\overline K_{\0}, \overline K_{\gk}$ and the functions $\frc_{\0},\frc_{\gk},\frs_{\gk}$ for $\gk\in\Z_{0}^{d}$ are respectively the eigenvalues and the eigenvectors of the operator on $\H$ defined by \eref{eq-oper}. 

Under our Assumption~\ref{H4}, this kernel is of the form \eref{eq-Kgene} and fits in our general framework. We may take $\frJ=\{\0\}\cup (\Z_{0}^{d}\times\{0,1\})$, $\fru_{\0}=1$, $l_{\0}=\overline K_{\0}=\sin^{2}s_{\0}\in [0,1]$  and for every $\gk\in \Z_{0}^{d}$  $\fru_{\gk,0}=\frc_{\gk}$, $\fru_{\gk,1}=\frs_{\gk}$ and $s_{\gk,0}=s_{\gk,1}$ satisfy $\overline K_{\gk}=\sin^{2}s_{\gk,0}=l_{\gk,0}=l_{\gk,1}\in [0,1]$. The condition 
\[
\sum_{j\in \frJ}l_{j}=\overline K_{\0}+2\sum_{\gk\in\Z_{0}^{d}}\overline K_{\gk}<+\infty
\]
is satisfied under Assumption~\ref{H2}. 

We may therefore  associate $\overline K$ with a DPP $X$ whose distribution  is denoted by $P_{\overline K}$. By \eref{eq-nb-pt}, it satisfies 
\begin{equation}\label{eq-ptmoy}
\E\cro{X(\T^{d})}=\int_{\T^{d}}K(x,x)d\mu(x)=\overline K(0)>0
\end{equation}
which implies that $X$ is not a.s.\ empty. 

We may alternatively write the distribution $P_{\overline K}$ as  $P_{\gs,\fru}$ (see \eref{def-Psu}) where $\fru=(u_{j})_{j\in\frJ}$ is known and given by the trigonometric functions defined above whereas $\gs\in[0,\pi/2]^{\frJ}$ is actually the only unknown parameter. 

If $I$ is a finite subset of $\Z_{0}^{d}$, the function 
\[
\overline K_{I}(x)=\overline K_{\0}+2\sum_{\gk\in I}\overline K_{\gk}\cos(2\pi \scal{\gk}{x})\quad \text{for $x\in \T^{d}$,}
\]
also belongs to $\cK$ and is associated with the distribution $P_{\overline K_{I}}=P_{\gs_{J},\fru_{J}}$ where  $J=\{\0\}\cup (I\times \{0,1\})\subset \frJ$ and $\gs_{J}=(s_{j}\1_{j\in J})_{j\in\frJ}$. Using \eref{eq-prop-h-DPP0} and the inequality $0\le x\le (\pi/2)\sin x$ for $x\in [0,\pi/2]$, we know that $h^{2}(P_{\overline K},P_{\overline K_{I}})=h^{2}(P_{\gs,\fru},P_{\gs_{J},\fru_{J}})$ is not larger than 
\begin{align}\label{eq-approx00}
\frac{1}{2}\sum_{j\not \in J}s_{j}^{2}=\sum_{\gk\not \in I}s_{\gk,0}^{2}\le \frac{\pi^{2}}{4}\sum_{\gk\not \in I}\sin^{2}s_{\gk,0}= \frac{\pi^{2}}{4}\sum_{\gk\not \in I}\overline K_{\gk}= \frac{\pi^{2}}{8}\pa{\overline K-\overline K_{I}}(0).
\end{align}
This result shows that the distribution $P_{\overline K_{I}}$ is close to $P_{\overline K}$ as soon as $\overline K_{I}$ approximates $\overline K$ sufficiently well at 0. To quantify this approximation error, we introduce the modulus of continuity of $\overline K$ at 0 defined by 
\begin{equation}\label{def-w}
w(z)=\sup_{|x|_{\infty}\le z}\ab{\overline K(x)-\overline K(0)}\quad \text{for $z\ge 0$},
\end{equation}
where $|x|_{\infty}=\max_{i\in\{1,\ldots,d\}}|x_{i}|$ for $x\in\T^{d}$. We say that $\overline K$ is $(L,\alpha)$-H\"older continuous at 0 if $w(z)\le Lz^{\alpha}$ for some $L,\alpha>0$. 
Note that if $\overline K$ is differentiable at 0, its gradient  $\nabla_{0}\overline K $ at $0\in\T^d$ is necessarily 0 since $\overline K$ is even. Consequently, if the gradient $x\mapsto \nabla_{x}\overline K$ of $\overline K$ exists and is $(L,\alpha')$-H\"older continuous at 0, in the sense that 
\begin{equation}\label{eq-Holder}
\sup_{|x|_{\infty}\le z}\ab{\nabla_{x}\overline K-\nabla_{0}\overline K}_{1}\le Lz^{\alpha'}\quad \text{for $L>0$ and $\alpha'\in (0,1]$,}
\end{equation}
we automatically get 
\begin{align*}
w(z)&=\sup_{|x|_{\infty}\le z}\ab{\overline K(x)-\overline K(0)-\scal{\nabla_{0}\overline K}{x}}=\sup_{|x|_{\infty}\le z}\ab{\int_{0}^{1}\scal{\nabla_{tx}\overline K-\nabla_{0}\overline K}{x}dt}\\
&\le z\sup_{|x|_{\infty}\le z}\ab{\nabla_{x}\overline K-\nabla_{0}\overline K}_{1}\le Lz^{1+\alpha'}.
\end{align*}
The inequality $w(z)\le Lz^{\alpha}$  is therefore not restricted to values of $\alpha\in (0,1]$.

For $D\in\N\et$, let us set 
\begin{equation}\label{def-scD}
\sC(D)=\{\gk\in\Z_{0}^d,\; |\gk|_{\infty}\le D\} \quad \text{so that}\quad |\sC(D)|= \frac{(2D+1)^{d}-1}{2}.
\end{equation}
The following result is proven in Section~\ref{pf-prop-approx}.
\begin{prop}\label{prop-approx}
For $D\ge 1$ and $I=\sC(D)$, 
\begin{align}\label{eq-approx}
\pa{\overline K-\overline K_{I}}(0)\le \inf_{\delta>(2\pi D)^{-1}}\frac{w(\delta)}{1-(2\pi D\delta)^{-1}}\le 2w\pa{\frac{1}{\pi D}}.
\end{align}
In particular, if $w(z)\le Lz^{\alpha}$ for some positive numbers $L$ and $\alpha$, 
\begin{equation}\label{eq-approx01}
\pa{\overline K-\overline K_{I}}(0)\le  C(\alpha)LD^{-\alpha}\quad \text{with}\quad C(\alpha)= \frac{(1+\alpha)^{1+\alpha}}{(2\pi \alpha)^{\alpha}}.
\end{equation}
\end{prop}

\subsection{Estimating the distribution $P_{\overline K}$ for $\overline K\in\cK$}
As already mentioned, the distribution of $P_{\overline K}$ is of the form $P_{\gs,\fru}$ where $\fru\in\U$ is known and we only need to estimate the parameter $\gs$ which is of the form $(s_{\0},(s_{\gk,0},s_{\gk,0})_{\gk\in \Z_{0}^{d}})$ with $s_{\0}$ and $s_{\gk,0}$ in $[0,\pi/2]$ for every $\gk\in\Z_{0}^{d}$. To do so, we apply our Corollary~\ref{cor-su} with $J=J(D)$ of the form 
\begin{equation}\label{def-JD}
J=J(D)=\{\0\}\cup\pa{\sC(D)\times\{0,1\}}\in\cJ\quad \text{for $D\ge 1$.}
\end{equation}
By \eref{eq-approx00} and \eref{eq-approx}, the parameter $\gt=\gs_{J(D)}\in \gS_{J(D)}$ satisfies 
\[
h^{2}\pa{P_{\gs,\fru},P_{\gs_{J(D)},\fru_{J(D)}}}\le \frac{1}{2}\sum_{j\not \in J(D)}s_{j}^{2}\le \frac{\pi^{2}}{4}w\pa{\frac{1}{\pi D}}.
\]
Identifying $\gS_{J(D)}$ with the set 
\[
[0,\pi/2]\times ([0,\pi/2]\times [0,\pi/2])^{\sC(D)}
\]
we may equip it with the prior given by the image $\gpi_{1,J(D)}$ of the uniform distribution on $[0,\pi/2]\times [0,\pi/2]^{\sC(D)}$ by the mapping $\gt=(t_{\0},(t_{\gk})_{\gk\in\sC(D)})\mapsto (t_{\0},(t_{\gk},t_{\gk})_{\gk\in\sC(D)})$. 

Since $\fru$ is known, for $J$ given by \eref{def-JD}, we take for $E$ the linear span of $\frc_{\0},\fru_{\gk,0}=\frc_{\gk}$ and $\fru_{\gk,1}=\frs_{\gk}$ for $\gk\in\sC(D)$ and equip $\U_{J}(E)$ with the Dirac mass $\delta_{\fru_{J(D)}}$ at $\fru_{J(D)}$. 

We then derive from Corollary~\ref{cor-su} the following result. The details of the proof are given in Section~\ref{pf-prop-nonadapt}.
\begin{prop}\label{prop-nonadapt}
Let $D\in \N\et$, $\overline K\in\cK$ with modulus of continuity $w$ at 0 given by \eref{def-w} and $\etc{X}$ be an $n$-sample with distribution $P_{\overline K}$ with $n\ge 2$. If $\gpi$ is the prior on $\Theta=\gS\times\U$ defined by $\gpi=\gpi_{1,J(D)}\otimes \delta_{\fru_{J(D)}}$, then for every $\xi>0$, with probability at least $1-2e^{-\xi}$, a $T$-Bayes estimator $(\widehat \gs,\fru_{J(D)})$ with distribution $\gpi_{\bsX}$ satisfies 
\begin{align*}
C^{-1}h^{2}\pa{P_{\gs,\fru},P_{\smash{\widehat \gs,\fru_{J(D)}}}}&\le w\pa{\frac{1}{\pi D}}+\frac{(1+2D)^{d}}{n}\log n+\frac{\xi}{n}.
\end{align*}
In particular, 
\begin{equation}\label{es-riskKbar}
C^{-1}\E\cro{h^{2}\pa{P_{\gs,\fru},P_{\smash{\widehat \gs,\fru_{J(D)}}}}}\le w\pa{\frac{1}{\pi D}}+\frac{(1+2D)^{d}}{n}\log n.
\end{equation}
\end{prop}
We can establish a bound which is adaptive with respect to the choice of $D$ by considering the hierarchical prior 
\begin{equation}\label{def-piKbar}
\gpi=\frac{6}{\pi^{2}}\sum_{D=1}^{+\infty}\frac{1}{D^{2}}\cro{\gpi_{1,J(D)}\otimes \delta_{\fru_{J(D)}}}
\end{equation}
 By taking $\cJ'=\{J(D),\; D\ge 1\}$ and $\cM_{J(D)}=\{D\}$ for every $D\ge 1$, 
we obtain that $\gpi$ is of the form \eref{def-hprior} with $\nu(J(D),D)=6/(\pi^{2} D^{2})$. Corollary~\ref{cor-su01} therefore applies. Since the complexity term $-\log \nu(J(D),D)=\log(\pi^{2}/6)+2\log D$ is negligible compared with $(1+2D)^{d}$, we derive from this corollary the following result.
\begin{prop}\label{prop-adapt}
Let $\gpi$ be the prior defined by \eref{def-piKbar}, $n\ge 2$, and $\etc{X}$ be an $n$-sample with distribution $P_{\overline K}$ for $\overline K\in\cK$.  For every $\xi>0$, with probability at least $1-2e^{-\xi}$, a $T$-Bayes estimator $(\widehat \gs,\fru_{J(\widehat D)})$ with distribution $\gpi_{\bsX}$ satisfies 
\begin{align}
h^{2}\pa{P_{\gs,\fru},P_{\smash{\widehat \gs,\fru_{J(\widehat D)}}}}&\le C\pa{B_{n}+\frac{\xi}{n}}\quad \text{with}\quad B_{n}=\inf_{D\ge 1}\cro{w\pa{\frac{1}{\pi D}}+\frac{(1+2D)^{d}}{n}\log n}.\label{def-Bn}
\end{align}
where $w$  is the modulus of continuity of $\overline K$ at 0 given by \eref{def-w}. In particular, 
\begin{align*}
\E\cro{h^{2}\pa{P_{\gs,\fru},P_{\smash{\widehat \gs,\fru_{J(\widehat D)}}}}}&\le CB_{n}.
\end{align*}
\end{prop}
If $\overline K$ is $(L,\alpha)$-H\"older continuous at 0, the risk of our randomised estimator is not larger than 
\begin{equation}\label{borne-Bn}
C(d)\cro{L^{d/(\alpha+d)}\pa{\frac{\log n}{n}}^{\alpha/(\alpha+d)}+\frac{\log n}{n}}
\end{equation}
for some positive constant $C(d)$ only depending on $d$. We recall that if $\overline K$ is differentiable on $\T^{d}$ with a $(L,\alpha')$-H\"older gradient at 0, the above inequality applies with $\alpha=1+\alpha'$. 

\subsection{Estimating $\overline K$}
In this section, we tackle the problem of identifying and estimating the parameter $\overline K$. In the sequel, we denote by $\cK_{+}$ the subset of $\cK$ consisting of the elements $\overline K$ which are nonnegative on $\T^d$. The following result is proven in Section~\ref{pf-prop-minh}. 
\begin{prop}\label{prop-minh}
Let $\overline K$ and $\overline L$ be two elements of $\cK$. We set
\begin{align*}
M(0)=\pa{\overline K(0)-\|\overline K\|_{2}^{2}}\vee \pa{\overline L(0)-\|\overline L\|_{2}^{2}}\quad \text{and}\quad \overline M(0)=\overline K(0)\vee \overline L(0)\ge M(0).
\end{align*}
Then, when $M(0)>0$ 
\begin{equation}\label{eq-prop-minh}
h^{2}\pa{P_{\overline K},P_{\overline L}}\ge 
 \cro{ \frac{\pa{\overline K(0)-\overline L(0)}^{2}}{10M(0)}\vee \frac{\|\overline K^{2}-\overline L^2\|_{2}^{2}}{191\overline M^{2}(0)\pa{\overline M^{2}(0)\vee 1}}}\wedge \frac{1}{10}.
\end{equation}
When $M(0)=0$, $P_{\overline K}$ and $P_{\overline L}$ are the distributions of two determinantal projection processes and $h^{2}\pa{P_{\overline K},P_{\overline L}}\ge 1/2$ when $\overline K(0)\ne \overline L(0)$ while 
\begin{equation}\label{eq-prop-minh01}
h^{2}\pa{P_{\overline K},P_{\overline L}} \ge \frac{\|\overline K^{2}-\overline L^2\|_{2}^{2}}{56\overline M^{2}(0)\pa{\overline M^{2}(0)\vee 1}}\quad \text{otherwise}.
\end{equation}
In particular, the parameter $\overline K$ is identifiable within $\cK_{+}$. 
\end{prop}
Let us now comment on these results. Inequality \eref{eq-prop-minh} implies that for $M(0)>0$, if $P_{\overline K}$ and $P_{\overline L}$ are close, so are $\overline K(0)$ and $\overline L(0)$. By~\eref{eq-ptmoy}, this means that the average number of points of a DPP with distribution $P_{\overline K}$ is close to that of a DPP with  distribution $P_{\overline L}$. When $M(0)=0$, it follows from \eref{eq-coefF} that $\overline K_{\gk}$ and $\overline L_{\gk}$ belong to $\{0,1\}$ for every $\gk\in \Z^{d}$. By \eref{eq-KFour}, this means that $P_{\overline K}$ and $P_{\overline L}$ are the distributions of two determinantal projection processes and their cardinalities are therefore deterministic.

When $M(0)>0$, the inequality 
\begin{equation}\label{eq-lowb01}
h^{2}\pa{P_{\overline K},P_{\overline L}}\ge \frac{1}{10}\cro{\frac{\pa{\overline K(0)-\overline L(0)}^{2}}{M(0)}\wedge 1},
\end{equation}
which is a consequence of \eref{eq-prop-minh}, is nearly sharp. Let us consider the case where $\overline K=\sin^{2} s$ and $\overline L=\sin^{2} (s+\epsilon)$ are constant on the torus $\T^{d}$ with $s\in (0,\pi/2)$, $\epsilon>0$ tending to 0. On the one hand, using that $\sin^{2}\! (s+\epsilon)-\sin ^{2}s=(2\cos s\sin s)\epsilon(1+o(1))$, we obtain that 
\begin{align*}
\frac{\pa{\overline K(0)-\overline L(0)}^{2}}{M(0)}&=\frac{\ab{\sin^{2} s-\sin^{2}(s+\epsilon)}^{2}}{\max\{\sin^{2} s-\sin^{4} s,\sin^{2}(s+\epsilon) -\sin^{4} (s+\epsilon)\}}\\
&=  \frac{(4\sin^{2}s\cos^{2}s)\epsilon^{2}}{\sin^{2} s\cos^{2}s}(1+o(1))=4\epsilon^{2}(1+o(1)).
\end{align*}
On the other hand, denoting by $\cU(\T^{d})$ the uniform distribution on $\T^{d}$, it follows from \eref{def-Psu} that 
\[
P_{\overline K}=\cos^{2}s\; \delta_{\vide}+\sin^{2}s\; \cU(\T^{d})\quad \text{and}\quad P_{\overline L}=\cos^{2} (s+\epsilon)\;\delta_{\vide}+\sin^{2}(s+\epsilon)\; \cU(\T^{d})
\]
so that
\begin{align*}
h^{2}(P_{\overline K},P_{\overline L})&=1-\cro{\cos s\cos(s+\epsilon)+\sin s\sin (s+\epsilon)}=1-\cos \epsilon=\frac{\epsilon^{2}}{2}(1+o(1)).
\end{align*}
We conclude that the ratio
\[
\frac{h^{2}(P_{\overline K},P_{\overline L})}{\pa{\overline K(0)-\overline L(0)}^{2}/M(0)}=\frac{1+o(1)}{8},
\]
and the lower bound \eref{eq-lowb01} that we get only differs by the factor $8/10=4/5$.

Combining the results of Propositions~\ref{prop-adapt} and~\ref{prop-minh} allows us to provide an estimator of $\overline K$ and a risk bound for this estimator by arguing as follows. Equipping the parameter space  $\gS\times\U$ with the hierarchical prior $\gpi$ defined by \eref{def-piKbar}, our $T$-Bayes estimator returns both an integer $\widehat D\ge 1$ and an element $\widehat \gs=(\widehat s_{\0},(\widehat s_{\gk},\widehat s_{\gk})_{\gk\in \sC(\widehat D)})\in\gS$. We may then define our estimator of $\overline K$ as the function
\begin{equation}\label{def-Khat}
\widehat K(x)=\widehat K_{\0}+2\sum_{\gk\in\sC(\widehat D)}\widehat K_{\gk}\cos\pa{2\pi\scal{\gk}{x}}\quad \text{for $x\in\T^{d}$}
\end{equation}
with $\widehat K_{\0}=\sin^{2}\widehat s_{\0}$ and $\widehat K_{\gk}=\sin^{2}\widehat s_{\gk}$ for every $\gk\in\sC(\widehat D)$. The function $\widehat K$ belongs to $\cK$ a.s.\ (since the distribution of $\widehat K_{\0}$ is absolutely continuous with respect to the Lebesgue measure) and $P_{\widehat K}=P_{\smash{\widehat \gs,\fru_{J(\widehat D)}}}$. We may then derive from Propositions~\ref{prop-adapt} and~\ref{prop-minh} the following result. The details of the proof are given in Section~\ref{pf-thm-borne}.

\begin{thm}\label{thm-borne}
Let $\gpi$ be the prior defined by \eref{def-piKbar}, $n\ge 2$, and $\etc{X}$ be an $n$-sample with distribution $P_{\overline K}$ for $\overline K\in\cK$.  The $T$-Bayes estimator $\widehat K$  
defined by \eref{def-Khat} possesses the following property. Given $\xi>0$, if $n$ is large enough to satisfy  the inequality $10C(B_{n}+\xi/n)< 1\wedge \overline K(0)$ where $C$ is the numerical constant in \eref{def-Bn}, with probability at least $1-2e^{-\xi}$, 
\begin{align}
\norm{\overline K^{2}-\widehat K^{2}}_{2}^{2}\le C'\cro{\overline K(0)\pa{1\vee \overline K(0)}}^{2}\cro{B_{n}+\frac{\xi}{n}}.\label{eq-thm-borne}
\end{align}
\end{thm}
The bound $B_{n}$ only depends on the smoothness of $\overline K$ at 0 and  we have seen that it is not larger than the quantity given in \eref{borne-Bn} when $\overline K$ is $(L,\alpha)$-H\"older continuous at 0. We note that the smaller $\overline K(0)=\E\cro{|X_{1}|}$, the smaller the risk bound in \eref{eq-thm-borne}. When $\overline K(0)\ge 1$,  the risk bound is proportional to $\overline K^{4}(0)$.

Of special interest is the case where $\overline K$ is positive on $\T^{d}$. Then  $\overline K$ is identifiable and $\fra=\min_{x\in\T^{d}}\overline K(x)$ is a positive number since $\overline K$ is continuous on $\T^{d}$. We deduce that 
\begin{align*}
\norm{\overline K^{2}-\widehat K^{2}}_{2}^{2}&=\int_{\T^{d}}\pa{\overline K^{2}-\widehat K^{2}}^{2}d\mu(x)=\int_{\T^{d}}\pa{\overline K-|\widehat K|}^{2}\pa{\overline K+|\widehat K|}^{2}d\mu(x)\\
&\ge \fra^{2}\int_{\T^{d}}\pa{\overline K-|\widehat K|}^{2}d\mu(x)=\fra^{2}\norm{\overline K-|\widehat K|}_{2}^{2}.
\end{align*}
As a consequence of \eref{eq-thm-borne}, the estimator $|\widehat K|$ of $\overline K$ satisfies 
\begin{align}
\norm{\overline K-|\widehat K|}_{2}^{2}\le C'\cro{\frac{\overline K(0)\pa{1\vee \overline K(0)}}{\fra}}^{2}\cro{B_{n}+\frac{\xi}{n}}\label{eq-thm-borne2}
\end{align}
and 
\begin{align}
\E\cro{\norm{\overline K-|\widehat K|}_{2}^{2}}\le C'\cro{\frac{\overline K(0)\pa{1\vee \overline K(0)}}{\fra}}^{2}B_{n}\label{eq-thm-borne2b}
\end{align}

and estimates $\overline K$ at rate $B_{n}$.

\section{Proofs}\label{sect-proof}
We start with some additional results about determinantal processes.  Given the kernel $K$, we may describe the distribution of a DPP $X$ as follows: for every $k\ge 1$ and all pairwise disjoint measurable subsets $A_{1},\ldots,A_{k}$ of $\cX$, 
\begin{align}\label{eq-def-PDP}
\E\cro{\prod_{i=1}^{k}X(A_{i})}=\int_{A_{1}\times\ldots\times A_{k}}\varrho_{k}(x_{1},\ldots,x_{k})d\mu(x_{1})\ldots d\mu(x_{k}).
\end{align}
where 
\[
0\le \varrho_{k}(x_{1},\ldots,x_{k})=\det\pa{K(x_{i},x_{j})}_{1\le i,j\le k}\quad \text{for every $x=(x_{1},\ldots,x_{k})\in\cX^{k}$.}
\]
In particular, the expected number of points in $X\cap A$ is 
\begin{equation}\label{eq-PDP-f0}
0\le \E\cro{X(A)}=\int_{A}K(x,x)d\mu(x).
\end{equation}
Another consequence of \eref{eq-def-PDP} is that for every nonnegative function $f$ on $\cX^k$ which vanishes whenever two of the coordinates of $x=(x_{1},\ldots,x_{k})$ coincide
\begin{equation}\label{eq-PDP-f}
\E\cro{\int_{\cX^{k}}f(x_{1},\ldots,x_{k})\prod_{i=1}^{k}dX(x_{i})}=\int_{\cX^k}f(x_{1},\ldots,x_{k})\varrho_{k}(x_{1},\ldots,x_{k})\prod_{i=1}^{k}d\mu(x_{i}).
\end{equation}
When $\cX$ is finite and equipped with the counting measure $\mu$, a consequence of \eref{eq-def-PDP} is that for every $A\subset \cX$ 
\[
\P\cro{A\subset X}=\E\cro{\prod_{a\in A}X(\{a\})}=\det\cro{(K(x,y))_{x,y\in A}}.
\]
We use the following results. 
\begin{prop}\label{prop-GS}
Let $v_{1},\ldots,v_{k}$ be $k\ge 1$ vectors lying on the unit sphere of a finite-dimensional linear subspace $E$ of $\H$. There  exists $\overline \gv=(\overline v_{1},\ldots,\overline v_{k})$ in the Stiefel manifold $\U_{k}(E)$ that satisfies 
\begin{equation}\label{eq-prop-GS}
\sum_{j=1}^{k}\norm{v_{j}-\overline v_{j}}^{2}=\inf_{\gu\in\U}\sum_{j=1}^{k}\norm{v_{j}-u_{j}}^{2}=\inf_{\gu\in\U_{k}(E)}\sum_{j=1}^{k}\norm{v_{j}-u_{j}}^{2}
\end{equation}
and this element $\overline \gv$ is unique when $v_{1},\ldots,v_{k}$ are linearly independent. Besides, if $v_{1},\ldots,v_{k}$ are i.i.d.\ random vectors uniformly distributed on the unit sphere of $E$, $\overline \gv$  is distributed as the Haar measure $\frm_{k}^{E}$ on $\U_{k}(E)$. 
\end{prop}
\begin{proof}[Proof of Proposition~\ref{prop-GS}]
Let $D\ge k$ be the dimension of $E$. We first choose an orthonormal basis of $E$ and denote by  $V$ the $D\times k$ random matrix whose columns are the coordinates of $v_{1},\ldots,v_{k}$ in this basis. Using a singular value decomposition, the $D\times k$ matrix $V$ can be written in the form $V=W\Delta P$ where the columns of $W$ are orthonormal in $\R^{D}$, $\Delta$ is the diagonal matrix with diagonal entries $\sigma_{1},\ldots,\sigma_{k}$ and $P$ is a $k\times k$ orthogonal matrix. For $i\in\{1,\ldots,k\}$, let $w_{i}$ be the vector of $E$ which is represented (in the orthonormal basis of $E$) by the $i$-th column of $W$. Since $V=W\Delta P$, we may then write for every $j\in\{1,\ldots,k\}$, $v_{j}=\sum_{i=1}^{k}\sigma_{i}P_{i,j}w_{i}$ and we set $\overline v_{j}=\sum_{i=1}^{k}P_{i,j}w_{i}$. Since the matrix $P$ is orthogonal, the vectors $\overline v_{1},\ldots,\overline v_{k}$ are orthonormal in $E$ and they satisfy 
\begin{align}\label{eq-sigma}
\sum_{j=1}^{k}\scal{v_{j}}{\overline v_{j}}=\sum_{j=1}^{k}\cro{\sum_{i=1}^{k}\sigma_{i}P_{i,j}^{2}}=\sum_{i=1}^{k}\sigma_{i}\cro{\sum_{j=1}^{k}P_{i,j}^{2}}=\sum_{i=1}^{k}\sigma_{i}.
\end{align}

Let now $\gu=(u_{1},\ldots,u_{k},\ldots)\in \U$ be an orthonormal system of $\H$ and $\Pi$ the orthogonal projector onto $E$. Since $u_{1},\ldots,u_{k}$ are orthonormal, for every $i\in\{1,\ldots,k\}$, $\sum_{j=1}^{k}\scal{u_{j}}{w_{i}}^{2}\le \norm{w_{i}}^{2}=1$ and it follows thus from the Cauchy--Schwarz inequality that
\begin{align}
\sum_{j=1}^{k}\scal{v_{j}}{u_{j}}&=\sum_{j=1}^{k}\scal{v_{j}}{\Pi u_{j}}=\sum_{j=1}^{k}\cro{\sum_{i=1}^{k}\sigma_{i}P_{i,j}\scal{u_{j}}{w_{i}}}= \sum_{i=1}^{k}\sigma_{i} \sum_{j=1}^{k}P_{i,j}\scal{u_{j}}{w_{i}}\nonumber\\
&\le \sum_{i=1}^{k}\sigma_{i} \pa{\sum_{j=1}^{k}P_{i,j}^{2}}^{1/2}\pa{\sum_{j=1}^{k} \scal{u_{j}}{w_{i}}^{2}}^{1/2}\le  \sum_{i=1}^{k}\sigma_{i}=\sum_{j=1}^{k}\scal{v_{j}}{\overline v_{j}}.\label{eq-CS}
\end{align}
We deduce that 
\begin{align*}
\frac{1}{2}\sum_{j=1}^{k}\norm{v_{j}-u_{j}}^{2}=k-\sum_{j=1}^{k}\scal{v_{j}}{u_{j}}\ge k-\sum_{j=1}^{k}\scal{v_{j}}{\overline v_{j}}=\frac{1}{2}\sum_{j=1}^{k}\norm{v_{j}-\overline v_{j}}^{2},
\end{align*}
which yields \eref{eq-prop-GS}.

 Let us now prove uniqueness when $v_{1},\ldots,v_{k}$ are linearly independent. In this case, the singular values $\sigma_{1},\ldots,\sigma_{k}$ of $V$ are positive. Let $\gu=(u_{1},\ldots,u_{k})\in\U_{k}(E)$ such that $\sum_{j=1}^{k}\norm{v_{j}-u_{j}}^{2}=\sum_{j=1}^{k}\norm{v_{j}-\overline v_{j}}^{2}$. Since the $\sigma_{i}$ are positive, it follows from \eref{eq-CS} that this equality is possible if and only if there is an equality when applying the Cauchy--Schwarz inequality and when $\sum_{j=1}^{k}\scal{u_{j}}{w_{i}}^{2}=1$ for every $i\in\{1,\ldots,k\}$. If there is equality when applying the Cauchy--Schwarz inequality, for every $i\in\{1,\ldots,k\}$, $P_{i,j}=c_{i}\scal{u_{j}}{w_{i}}$ or $c_{i}P_{i,j}=\scal{u_{j}}{w_{i}}$ for some $c_{i}\in\R_{+}$  for every $j\in\{1,\ldots,k\}$. Since $\sum_{j=1}^{k}P_{i,j}^{2}=1=\sum_{j=1}^{k}\scal{u_{j}}{w_{i}}^{2}$, we necessarily have  $c_{i}=1$ and $P_{i,j}=\scal{u_{j}}{w_{i}}$ for every $i,j\in\{1,\ldots,k\}$. Besides, the identity $\norm{w_{i}}^{2}=1=\sum_{j=1}^{k}\scal{u_{j}}{w_{i}}^{2}$ implies that $w_{i}$ belongs to the linear span $U$ of $u_{1},\ldots,u_{k}$ for  every $i\in\{1,\ldots,k\}$. Hence $w_{1},\ldots,w_{k}$ is also an orthonormal basis of $U$. Consequently, $\overline v_{j}=\sum_{i=1}^{k}P_{i,j}w_{i}=\sum_{i=1}^{k}\scal{u_{j}}{w_{i}}w_{i}=u_{j}$ for every $j\in\{1,\ldots,k\}$. 
 
Let us now prove that $\overline \gv$ has distribution $\frm_{k}^{E}$ when $v_{1},\ldots,v_{k}$ are i.i.d.\ and uniformly distributed on the unit sphere of $E$. Note that the $v_{j}$ are almost surely linearly independent when $k\le D$. We recall that $\frm_{k}^{E}$ is the unique probability measure on $\U_{k}(E)$ which is invariant under the mapping $\gu\in \U_{k}(E)\mapsto O\gu=(Ou_{1},\ldots,Ou_{k})$ for every orthogonal transformation $O$ of $E$. It is therefore sufficient to prove that $\overline \gv$ and $O\overline \gv=(O\overline v_{1},\ldots,O\overline v_{k})$ have the same distribution for every orthogonal transformation $O$. 

Given $O$, let us observe that 
\begin{align*}
\inf_{\gu\in\U_{k}(E)}\sum_{j=1}^{k}\norm{Ov_{j}-u_{j}}^{2}&=\inf_{\gu\in\U_{k}(E)}\sum_{j=1}^{k}\norm{Ov_{j}-Ou_{j}}^{2}=\inf_{\gu\in\U_{k}(E)}\sum_{j=1}^{k}\norm{v_{j}-u_{j}}^{2}\\
&= \sum_{j=1}^{k}\norm{v_{j}-\overline v_{j}}^{2}=\sum_{j=1}^{k}\norm{Ov_{j}-O\overline v_{j}}^{2}.
\end{align*}
Consequently, $\argmin_{\gu\in\U_{k}(E)}\sum_{j=1}^{k}\norm{Ov_{j}-u_{j}}^{2}=O\overline \gv$. Since $v_{1},\ldots,v_{k}$ are i.i.d.\ and uniformly distributed on the unit sphere of $E$, $\gv=(v_{1},\ldots,v_{k})$ and $O\gv$ have the same distribution and so have the two random variables $\argmin_{\gu\in\U_{k}(E)}\sum_{j=1}^{k}\norm{v_{j}-u_{j}}^{2}=\overline \gv$ and $\argmin_{\gu\in\U_{k}(E)}\sum_{j=1}^{k}\norm{Ov_{j}-u_{j}}^{2}=O\overline \gv$.
\end{proof}

\begin{lem}\label{res-ine}
Let $a$ and $b$ be positive numbers and $c\in\R$. The inequality $ax+b\log x+c\ge 0$ is satisfied if 
\[
x\ge \frac{b}{a}\log_{+}\cro{\frac{a}{b}e^{-c/b}}.
\]
\end{lem}

\begin{proof}[Proof of Lemma~\ref{res-ine}]
Let us first note that $ax+b\log x+c\ge 0$ is equivalent to $(a/b)x\ge \log(e^{-c/b}/x)$ since $b>0$. Setting $y=e^{-c/b}/x$ and $u=(a/b)e^{-c/b}>0$, the last inequality is equivalent to $y\log y\le u$. Let $f(y)=y\log y$ for $y>0$ and $y\et$ be the solution of $f(y)=u$. Since $u>0$, $y\et>1$ and since $f$ is increasing on $[1,+\infty)$ and negative on $(0,1)$, $f(y)\le u$ is equivalent to $y\le y\et$. The function $f$ is convex on $(0,+\infty)$ and satisfies $f(1)=0$ and $f(e)=e$. Hence, $f(y)\le e(y-1)/(e-1)$ for every $y\in [1,e]$. In particular, if $u\in [0,e]$ and $y\le 1+(e-1)u/e$, $f(y)\le u$. Let us now assume that $u>e$. We note that $f(y\et)=u$ implies that $y\et\ge e$, hence, $\log u=\log f(y\et)=\log y\et+\log\log y\et\ge \log y\et$. We deduce that $u=y\et\log y\et\le y\et\log u$, that is, $y\et\ge u/\log u>0$. Hence,  for $u>e$, if $y\le u/\log u=u/\log_{+}u$, $f(y)\le f(y\et)= u$.  Since for all $u\in (0,e]$, $u/\log_{+}u =u\le  1+(e-1)u/e$, we conclude that $f(y)\le u$ if $y\le u/\log_{+}u$, which yields the conclusion. 
\end{proof}

\begin{lem}\label{lem-h}
Let $X,Y$ be two random variables with distributions $P$ and $Q$ respectively on a measurable space $E$ and $T$ be a measurable function from $E$ to $\R$ such that $\E\cro{T(X)}\ne \E\cro{T(Y)}$. Then, 
\begin{align}\label{eq-lem-h0}
h^{2}(P,Q)\ge \frac{\Delta^2}{4\Var\cro{T(X)}+4\Var\cro{T(Y)}+2\Delta^{2}}\quad \text{with}\quad \Delta=\E\cro{T(X)}-\E\cro{T(Y)}.
\end{align}
In particular, 
\begin{align}\label{eq-lem-h1}
h^{2}(P,Q)\ge \frac{\Delta^2}{4\E\cro{T^{2}(X)}+4\E\cro{T^{2}(Y)}}.
\end{align}
\end{lem}

\begin{proof}[Proof of Lemma~\ref{lem-h}]
Let $\mu$ be a measure that dominates both $P$ and $Q$ and let $p$ and $q$ be their respective densities. Then, by the Cauchy--Schwarz inequality, 
\begin{align*}
\ab{\E\cro{T(X)}-\E\cro{T(Y)}}^{2}&=\ab{\int_{E}T(x)\pa{p(x)-q(x)}d\mu(x)}^{2}\\
&=\ab{\int_{E}T(x)\pa{\sqrt{p(x)}+\sqrt{q(x)}}\pa{\sqrt{p(x)}-\sqrt{q(x)}}d\mu(x)}^{2}\\
&\le \int_{E}T^{2}(x)\pa{\sqrt{p(x)}+\sqrt{q(x)}}^{2}d\mu(x)\int_{E}\pa{\sqrt{p(x)}-\sqrt{q(x)}}^{2}d\mu(x)\\
&\le 4\cro{\int_{E}T^{2}(x)\pa{p(x)+q(x)}d\mu(x)}h^{2}(P,Q).
\end{align*}
It remains to apply this inequality to $T-c$ and to optimise with respect to $c$. The right-hand side is minimum at $c=\pa{\E[T(X)]+\E[T(Y)]}/2$.
\end{proof}

\subsection{Proof of Proposition~\ref{prop-Prop}}\label{pf-prop-Prop}
For \ref{Prop-i} and \ref{Prop-ii}, we refer to Baraud and Romon~\citeyearpar{Baraud:2026aa}. Let  $\theta\in \prod_{i=1}^{k}\Theta_{i}$. It follows from \eref{eq-hleHi} that if $r\ge \sum_{i=1}^{k}R_{i}(\gpi_{i},\theta_{i})$
\begin{align*}
\sB(\theta,r)\supset \prod_{i=1}^{k}\ac{\theta_{i}'\in \Theta_{i}, \gell_{i}(\theta_{i},\theta_{i}')\le R_{i}(\gpi_{i},\theta_{i})}.
\end{align*}
Consequently, by definition of the numbers $R_{i}(\gpi_{i},\theta_{i})$, for $r\ge \sum_{i=1}^{k}R_{i}(\gpi_{i},\theta_{i})$,
\begin{align*}
\gpi\pa{\sB(\theta,r)}&\ge \prod_{i=1}^{k}\gpi_{i}\pa{\ac{\theta_{i}'\in \Theta_{i}, \gell_{i}(\theta_{i},\theta_{i}')\le R_{i}(\gpi_{i},\theta_{i})}}\\
&\ge \prod_{i=1}^{k}\exp\pa{-\gamma  R_{i}(\gpi_{i},\theta_{i})}\ge \exp\pa{-\gamma r},
\end{align*}
which yields \eref{eq-hleHic}. 

\subsection{Proof of Proposition~\ref{prop-DPP}}\label{pf-prop-DPP}
We only need to prove the triangle inequality since the other properties are clearly satisfied. Let us consider $U,V,W$ three linear subspaces of $\H$ of dimension $k\ge 1$. We may consider a finite-dimensional linear space $E$ generated by an orthonormal basis $e_{1},\ldots,e_{D}$, $D\ge k$, that contains the three of them. Given $u_{1},\ldots,u_{k}\in E$ we denote by $\wedge \gu=u_{1}\wedge \ldots \wedge u_{k}$ the $k$-linear alternating form defined by
\[
\wedge\gu(v_{1},\ldots,v_{k})=\det \pa{G(\gu,\gv)}\quad \text{with}\quad G(\gu,\gv)=\pa{\scal{u_{j}}{v_{j'}}}_{1\le j,j'\le k}.
\]
Let $\mathop{\bigwedge}^{k}(E)$ be the linear span of the $k$-linear alternating forms $\wedge\gu$ for $u_{1},\ldots,u_{k}\in E$. Given a subset $S=\{i_{1},\ldots,i_{k}\}\subset \{1,\ldots,D\}$ with $1\le i_{1}<\ldots<i_{k}\le D$ we set
\[
\mathop{\bigwedge}_{i\in S}u_{i}=u_{i_{1}}\wedge\ldots\wedge u_{i_{k}}\quad \text{for $u_{1},\ldots,u_{D}\in E$.}
\]
Denoting by $I$ the set $\{1,\ldots,k\}$, it follows from the Cauchy--Binet formula (or simply the properties of the determinant) that 
\begin{align*}
\wedge\gu=u_{1}\wedge \ldots \wedge u_{k}=\sum_{S}\det\pa{\pa{\scal{u_{i}}{e_{j}}_{(i,j)\in I \times S}}}\mathop{\bigwedge}_{i\in S}e_{i}
\end{align*}
where the sum runs over the elements $S$ of the class $\cS_{k}$ of the subsets of $\{1,\ldots,D\}$ with cardinality $k$. We obtain that the family $\{\mathop{\bigwedge}_{i\in S}e_{i},\; S\in\cS_{k}\}$ generates $\mathop{\bigwedge}^{k}(E)$. It is actually easy to see that it is a basis since for $1\le i_{1}<\ldots<i_{k}\le D$, $\pa{\mathop{\bigwedge}_{i\in S}e_{i}}(e_{i_{1}},\ldots,e_{i_{k}})=1$ if $S=\{i_{1},\ldots,i_{k}\}$ and 0 otherwise. We may then equip $\mathop{\bigwedge}^{k}(E)$ with an inner product $[\cdot,\cdot]$ for which $\{\mathop{\bigwedge}_{i\in S}e_{i},\; S\in\cS_{k}\}$ is an orthonormal basis.  
We deduce that if $\wedge\gu=u_{1}\wedge\ldots \wedge u_{k}$ and $\wedge\gv=v_{1}\wedge \ldots \wedge v_{k}$ are two vectors of $\mathop{\bigwedge}^{k}(E)$, 
\begin{align*}
\cro{\wedge\gu, \wedge\gv}&=\sum_{S}\det\pa{\pa{\scal{u_{i}}{e_{j}}_{(i,j)\in I\times S}}}\det\pa{\pa{\scal{v_{i}}{e_{j}}_{(i,j)\in I\times S}}}=\det \pa{G(\gu,\gv)}
\end{align*}
by the Cauchy--Binet formula. In particular if $u_{1},\ldots,u_{k}$ and $v_{1},\ldots,v_{k}$ are orthonormal bases of $U$ and $V$ respectively, $[\wedge\gu,\wedge\gu]=[\wedge\gv,\wedge\gv]=1$ and
\[
\cro{\wedge\gu-\wedge\gv,\wedge\gu-\wedge\gv}=\cro{\wedge\gu,\wedge\gu}+\cro{\wedge\gv,\wedge\gv}-2\cro{\wedge\gu,\wedge\gv}=2\pa{1-\det G(\gu,\gv)}.
\]
In particular, $\ab{\wedge \gu}_{\wedge}=\sqrt{[\wedge\gu,\wedge\gu]}$ defines a norm on $\mathop{\bigwedge}^{k}(E)$ and it satisfies the triangle inequality. Let us now consider three linear spaces $U,V,W$ of dimension $k$ generated by the orthonormal bases $\{u_{1},\ldots,u_{k}\}$, $\{v_{1},\ldots,v_{k}\}$ and $\{w_{1},\ldots,w_{k}\}$ respectively. Replacing $u_{1}$ into $-u_{1}$ and $v_{1}$ into $-v_{1}$ if necessary, we may assume that  $\det G(\gu,\gw)$ and $\det G(\gv,\gw)$ are nonnegative in which case, 
\begin{align*}
\frd(U,W)+\frd(W,V)&=\frac{1}{\sqrt{2}}\cro{\ab{\wedge \gu-\wedge \gw}_{\wedge}+\ab{\wedge \gw-\wedge \gv}_{\wedge}}\ge \frac{\ab{\wedge \gu-\wedge \gv}_{\wedge}}{\sqrt{2}}\\
&=\sqrt{1-\det G(\gu,\gv)}\ge \sqrt{1-|\det G(\gu,\gv)|}=\frd(U,V).
\end{align*}
This proves the triangle inequality.

Let us now prove \eref{eq-DPP-aff}. 
Let us set $G=\pa{\scal{u_{j}}{v_{j'}}}_{1\le j,j'\le k}$ and $0\le \fra_{1}\le \ldots\le \fra_{k}\le \pi/2$ the principal angles between $U$ and $V$. Since the eigenvalues of $GG^{t}$ are $\cos^{2}\fra_{1},\ldots,\cos^{2}\fra_{k}$, $\ab{\det G}=\sqrt{\det GG^{t}}=\prod_{i=1}^{k}\cos \fra_{i}$. Applying Proposition~3.4 in Baraud~\citeyearpar{Baraud:2012qy}, we obtain that $\rho\pa{P_{U},P_{V}}\ge \ab{\det G}$, hence 
\begin{align*}
h^{2}(P_{U},P_{V})=1-\rho(P_{U},P_{V})\le 1-\prod_{i=1}^{k}\cos \fra_{i}= \frd^{2}(U,V)
\end{align*}
which proves the first inequality in \eref{eq-DPP-aff}.

Let  $\{u_{1}\et,\ldots,u_{k}\et\}$, $\{v_{1}\et,\ldots,v_{k}\et\}$ be principal vectors of $U$ and $V$.  We may write $u_{i}=\sum_{j=1}^{k}\scal{u_{i}}{u_{j}\et}u_{j}\et$, $v_{i}=\sum_{j=1}^{k}\scal{v_{i}}{v_{j}\et}v_{j}\et$. Since $\{u_{1},\ldots,u_{k}\}$ and $\{v_{1},\ldots,v_{k}\}$ are orthonormal bases of $U$ and $V$ respectively, 
$\sum_{i=1}^{k}\scal{u_{i}}{u_{j}\et}^{2}=\norm{u_{j}\et}^{2}=1$ and 
\begin{align}\label{eq-DPP-001} 
\sum_{i=1}^{k}\scal{v_{i}}{v_{j}\et}^{2}= 1\quad \text{for every $j\in\{1,\ldots,k\}$.}
\end{align}
Using that  $\scal{u_{j}\et}{v_{j}\et}=\cos \fra_{j}$ for every $j\in\{1,\ldots,k\}$ and the Cauchy--Schwarz inequality, 
\begin{align}
\sum_{i=1}^{k}\scal{u_{i}}{v_{i}}&=\sum_{i=1}^{k}\sum_{j=1}^{k}\cos \fra_{j}\scal{u_{i}}{u_{j}\et}\scal{v_{i}}{v_{j}\et}=\sum_{j=1}^{k}\cos \fra_{j}\cro{\sum_{i=1}^{k}\scal{u_{i}}{u_{j}\et}\scal{v_{i}}{v_{j}\et}}\nonumber\\
&\le \sum_{j=1}^{k}\cos \fra_{j}\cro{\sqrt{\sum_{i=1}^{k}\scal{u_{i}}{u_{j}\et}^{2}}\sqrt{\sum_{i=1}^{k}\scal{v_{i}}{v_{j}\et}^{2}}}\le \sum_{j=1}^{k}\cos \fra_{j}.\label{eq-DPP-002} 
\end{align}
Hence, using the inequality $1-\prod_{i=1}^{k}a_{i}\le \sum_{i=1}^{k}(1-a_{i})$ which holds for $a_{1},\ldots,a_{k}$ in $[0,1]$, we get
\begin{align}
\frac{1}{2}\sum_{i=1}^{k}\norm{u_{i}-v_{i}}^{2}&=\sum_{i=1}^{k}\cro{1-\scal{u_{i}}{v_{i}}}\ge  \sum_{i=1}^{k}\cro{1-\cos \fra_{i}}=\frac{1}{2}\sum_{i=1}^{k}\norm{u_{i}\et-v_{i}\et}^{2} \label{eq-DPP-003} \\
&\ge 1-\prod_{i=1}^{k}\cos \fra_{i}=\frd^{2}(U,V).\nonumber
\end{align}
This proves the second inequality in \eref{eq-DPP-aff}.

\subsection{Proof of Proposition~\ref{prop-rUni}}\label{pf-prop-rUni}
Let us identify $E$ with $\R^{D}$ by choosing an orthonormal basis of $E$ and denote by $\nu_{k,D}$ the Haar measure on the Grassmannian $\cG_{k}(\R^{D})$. If $O$ is an orthogonal transformation of $\R^{D}$ and $OV$ denotes the image of $V$ by $O$, we observe that $\frd(U,OV)=\frd(O^{-1}U,V)$ where $O^{-1}U$ is the image of $U$ by $O^{-1}$. Since the probability $\nu_{k,D}$ is stable under any orthogonal transformation, we deduce that if $V$ has distribution  $\nu_{k,D}$, $\frd(U,OV)=\frd(O^{-1}U,V)$ has the same distribution as $\frd(U,V)$. For the study of the random variable $\frd(U,V)$ with $V\sim \nu_{k,D}$, we may therefore assume with no loss of generality that $U$ is the linear span of the $k$-first vectors of the canonical basis of $\R^{D}$, which we shall do hereafter.

Let us now consider the $(D\times k)$-matrix $Z$ whose entries are independent standard Gaussian random variables $\cN(0,1)$. We denote by $Z_{1}$ the $(k\times k)$-matrix that consists of the $k$ first rows of $Z$ while $Z_{2}$ consists of the $D-k$ remaining rows. We may decompose (in a unique way) $Z=HT$ where $H$ is a $D\times k$ matrix whose columns are orthonormal vectors of $\R^{D}$ and $T$ is a $k\times k$ upper-triangular matrix with positive diagonal elements. We have seen that the linear span $V$ of the $k$ columns of $H$ is distributed as $\nu_{k,D}$. Note that $Z_{1}=H_{1}T$ where $H_{1}$ is the $k\times k$ matrix that consists of the $k$ first rows of $H$. Since $U$ is the linear span of the $k$ first vectors of the canonical basis of $\R^{D}$, it follows from~\eref{eq-def-d0} and \eref{eq-def-d} that
\[
1-\frd^{2}(U,V)=|\det H_{1}|=\sqrt{\det\pa{H_{1}^{t}H_{1}}}. 
\]
We observe that $Z_{1}^{t}Z_{1}=T^{t}H_{1}^{t}H_{1}T$ and 
$Z^{t}Z=Z_{1}^{t}Z_{1}+Z_{2}^{t}Z_{2}=T^{t}H^{t}HT=T^{t}T$ since the columns of $H$ are orthonormal. 
Consequently, 
\[
\det\pa{Z_{1}^{t}Z_{1}}=\det\pa{T^{t}H_{1}^{t}H_{1}T}=\det\pa{T^{t}T}\det\pa{H_{1}^{t}H_{1}}
\]
which leads to 
\begin{align*}
\cro{1-\frd^{2}(U,V)}^{2}=\det\pa{H_{1}^{t}H_{1}}=\frac{\det\pa{Z_{1}^{t}Z_{1}}}{\det\pa{Z_{1}^{t}Z_{1}+Z_{2}^{t}Z_{2}}},
\end{align*}
where  $Z_{1}^{t}Z_{1}$ and $Z_{2}^{t}Z_{2}$ are independent random variables. By applying Theorem~10.5.3 and Corollary~10.5.4 in Muirhead~\citeyearpar{MR652932} with $n-p=m=k$ and $r=D-k$, we obtain that $\cro{1-\frd^{2}(U,V)}^{2}$ has the same distribution as $\prod_{i=1}^{a}B_{i}$ which leads to \eref{eq-rUni00}. In particular, 
\begin{align}\label{prod-bi}
\P\cro{\prod_{i=1}^{a}B_{i}\ge  \pa{1-t^{2}}^{2}}\ge \prod_{i=1}^{a}\P\cro{B_{i}\ge  \pa{1-t^{2}}^{2/a}}.
\end{align}
We recall that for each $i\in\{1,\ldots,a\}$, $B_{i}$ has  density $q_{i}(x)=\alpha C_{i}x^{\alpha_{i}-1}(1-x)^{\alpha-1}$ with 
$\alpha_{i}=i/2$, $\alpha=b/2$ and $C_{i}=\Gamma(\alpha_{i}+\alpha)/[\Gamma(\alpha_{i})\Gamma(\alpha+1)]$. For $i\in\{1,\ldots,a\}$, 
%
\begin{align}\label{eq-mino-2}
\P\cro{B_{i}\ge  \pa{1-t^{2}}^{2/a}}&=C_{i}I_{i}\quad \text{with}\quad I_{i}=\alpha
\int_{(1-t^{2})^{2/a}}^{1}x^{\alpha_{i}-1}(1-x)^{\alpha-1}dx.
\end{align}

Let us first consider the case where $\alpha\ge 1$. Doing the change of variables $x=(1-u)^{2/a}$, $dx=-(2/a)(1-u)^{2/a-1}du$ we obtain that for every $i\in\{1,\ldots,a\}$, 
\begin{align*}
I_{i}&=\frac{2\alpha}{a}\int_{0}^{t^{2}}\pa{1-u}^{2\alpha_{i}/a-1}\pa{1-(1-u)^{2/a}}^{\alpha-1}du\\
&=\frac{2\alpha}{a}\int_{0}^{t^{2}}\pa{1-u}^{2\alpha_{i}/a-1}\cro{\pa{1-(1-u)^{1/a}}\pa{1+(1-u)^{1/a}}}^{\alpha-1}du\\
&\ge \frac{2\alpha}{a}\int_{0}^{t^{2}}\pa{1-u}^{2\alpha_{i}/a-1}\pa{1-(1-u)^{1/a}}^{\alpha-1}du.
\end{align*}
We note that $2\alpha_{i}/a=i/a\le 1$ for every $i\in\{1,\ldots,a\}$ and $(1-u)^{1/a}\le 1-u/a$ for every $u\in (0,1)$. Consequently, for every $i\in\{1,\ldots,a\}$
\begin{align}\label{eq-mino-3}
I_{i}\ge \frac{2\alpha}{a}\int_{0}^{t^{2}}\pa{\frac{u}{a}}^{\alpha-1}du=2\pa{\frac{t}{\sqrt{a}}}^{2\alpha}.
\end{align}
Using the equalities $t\Gamma(t)=\Gamma(t+1)$ for $t>0$ and $\alpha_{i+2}=\alpha_{i}+1$ for $i\ge 1$,we observe that for every $i\in\{1,\ldots,a\}$, 
\begin{align*}
\frac{C_{i+2}}{C_{i}}=\frac{\Gamma(\alpha_{i}+\alpha+1)\Gamma(\alpha_{i})}{\Gamma(\alpha_{i}+1)\Gamma(\alpha_{i}+\alpha)}=1+\frac{\alpha}{\alpha_{i}}\ge 1+\frac{b}{i}\ge 2.
\end{align*}
Besides, $C_{2}=1$ and it follows from the inequality $\Gamma(x+1/2)\ge \Gamma(x+1)/\sqrt{x+1/2}$ for $x>0$ (see Wendel~\citeyearpar{MR29448} or Qi~\citeyearpar{Qi_2010} for a review of inequalities on the Gamma function)  that 
\[
C_{1}=\frac{\Gamma(\alpha+1/2)}{\Gamma(\alpha+1)\Gamma(1/2)}>\frac{1}{\sqrt{\pi(\alpha+1/2)}}=\sqrt{\frac{2}{\pi(b+1)}}.
\]
We deduce that $C_{2j-1}\ge 2^{j-1}C_{1}$ for $1\le j\le \floor{(a+1)/2}$ and $C_{2j}\ge 2^{j-1}$ for $1\le j\le \floor{a/2}$. It follows from \eref{eq-mino-2} and \eref{eq-mino-3} that when $a=2p$ is even, 
\begin{align*}
\prod_{i=1}^{a}\P\cro{B_{i}\ge  \pa{1-t^{2}}^{2/a}}&=\prod_{i=1}^{a}C_{i}I_{i}=\cro{2\pa{\frac{t}{\sqrt{a}}}^{2\alpha}}^{a}\prod_{i=1}^{p}C_{2i-1}\prod_{i=1}^{p}C_{2i}\\
&\ge \cro{2\pa{\frac{t}{\sqrt{a}}}^{2\alpha}}^{a}2^{2\sum_{i=1}^{p}(i-1)}C_{1}^{p}\\
&=\pa{\frac{t}{\sqrt{a}}}^{ab}2^{(a+2)a/4}\pa{\frac{2}{\pi(b+1)}}^{a/4}.
\end{align*}
If $a=2p-1$ is odd, we also obtain that 
\begin{align*}
\prod_{i=1}^{a}\P\cro{B_{i}\ge  \pa{1-t^{2}}^{2/a}}
&\ge \pa{\frac{t}{\sqrt{a}}}^{ab}2^{(a+1)^{2}/4}\pa{\frac{2}{\pi(b+1)}}^{(a+1)/4}.
\end{align*}
In any case, 
\begin{align*}
\prod_{i=1}^{a}\P\cro{B_{i}\ge  \pa{1-t^{2}}^{2/a}}\ge \pa{\frac{t}{\sqrt{a}}}^{ab}\pa{\frac{1}{\pi(b+1)}}^{(a+1)/4}.
\end{align*}
Consequently, 
\begin{align*}
&\frac{1}{ab}\log\cro{\prod_{i=1}^{a}\P\cro{B_{i}\ge  \pa{1-t^{2}}^{2/a}}}\ge \log \pa{\frac{t}{\sqrt{a}}}-\frac{(a+1)\log [\pi(b+1)]}{4ab}.
\end{align*}
Note that $b\ge 2$ (because $\alpha\ge 1$) and $a\ge 1$. Using that $u\mapsto u^{-1}\log[\pi(1+u)]$ is nonincreasing on $[2,+\infty)$, 
\begin{align*}
\frac{(a+1)\log [\pi(b+1)]}{4ab}\le \frac{(a+1)\log(3\pi)}{8a}\le \frac{\log(3\pi)}{4}.
\end{align*}
We deduce that 
\begin{align*}
\frac{1}{ab}\log\cro{\prod_{i=1}^{a}\P\cro{B_{i}\ge  \pa{1-t^{2}}^{2/a}}}&\ge \log \pa{\frac{t}{\sqrt{a}}}-\frac{\log(3\pi)}{4}\ge\log \pa{\frac{t}{2\sqrt{a}}}
\end{align*}
and establishes \eref{eq-rUni01} when $\alpha\ge 1$. Let us now consider the situation where $\alpha<1$. Then $\alpha=1/2$ ($D=2$ and $k=1$), $a=b=1$ and 
\begin{align*}
\P\cro{B_{1}\ge  \pa{1-t^{2}}^{2}}&=\frac{\Gamma(1)}{\Gamma^{2}(1/2)}
\int_{(1-t^{2})^{2}}^{1}x^{-1/2}(1-x)^{-1/2}dx=\frac{4\arcsin(t/\sqrt{2})}{\pi}.
\end{align*}
Using that $\arcsin u\ge u$ for every $u\in [0,1]$, \eref{eq-rUni01} follows. 

\subsection{Proof of Corollary~\ref{cor-cont01}}\label{pf-cor-cont01}
Let $U\in\cG_{k}(E)$. Given $r> 0$, it follows from~\eref{eq-DPP-aff} that the ball $\sB(U,r)=\{V\in \cG_{k}(\H),\; nh^{2}(P_{U},P_{V})\le r\}$ contains the set $\{V\in \cG_{k}(E),\; \frd^{2}(U,V)\le (r/n)\wedge 1\}$, hence 
\begin{align*}
\gpi\pa{\sB(U,r)}&\ge \gpi\pa{\ac{V\in \cG_{k}(E),\; \frd(U,V)\le 1\wedge \sqrt{r/n}}}.
\end{align*}
Let us first note that $-\log\gpi\pa{\sB(U,r)}=0$ for $r>n$ and consequently $-\log\gpi\pa{\sB(U,r)}\le \gamma r$ for $r>n$.
For $r\in (0,n]$, we may apply \eref{eq-rUni01} with $t=\sqrt{r/n}\in [0,1]$ and obtain 
\begin{align*}
-\log \gpi\pa{\sB(U,r)}&\le  \frac{k(D-k)}{2}\log_{+}\pa{\frac{4na'}{r}}\quad \text{with}\quad a'=k\wedge (D-k).
\end{align*}
If $D=k$, $-\log\gpi\pa{\sB(U,r)}\le 0\le \gamma r$ for all $r> 0$, hence $\overline r(\gpi,U)=0$. When $D>k$, by applying Lemma~\ref{res-ine} with $a=\gamma$, $b=k(D-k)/2>0$ and $c=-k(D-k)\log(4na')/2$ we obtain that $-\log \gpi\pa{\sB(U,r)}\le \gamma r$ as soon as 
\[
r\ge \frac{k(D-k)}{2\gamma}\log_{+}\pa{\frac{8\gamma n a'}{k(D-k)}}=\frac{k(D-k)}{2\gamma}\log_{+}\pa{\frac{8\gamma n }{k\vee (D-k)}},
\]
which concludes the proof.

\subsection{Proof of Proposition~\ref{prop-approx01}}\label{pf-prop-approx01}
For $i\in\{1,\ldots,k\}$, let $v_{i}$ be the orthogonal projection of $u_{i}$ onto $E$, $V$ a $k$-dimensional linear subspace of $E$ that contains all the $v_{i}$ and $\Pi_{U}$ the orthogonal projection onto $U$. Let $\{u_{1}\et,\ldots,u_{k}\et\}$, $\{v_{1}\et,\ldots,v_{k}\et\}$ be principal vectors of $U$ and $V$ and $\fra_{1},\ldots,\fra_{k}$ the principal angles. With our definition of $v_{i}$, $v_{1},\ldots,v_{k}$ may not be orthonormal and \eref{eq-DPP-001} becomes  
\begin{align*}
\sum_{i=1}^{k}\scal{v_{i}}{v_{j}\et}^{2}=\sum_{i=1}^{k}\scal{u_{i}}{v_{j}\et}^{2}=\sum_{i=1}^{k}\scal{u_{i}}{\Pi_{U}v_{j}\et}^{2}=\norm{\Pi_{U}v_{j}\et}^{2}\le 1
\end{align*}
for every $j\in\{1,\ldots,k\}$. Arguing as for the proof of inequality \eref{eq-DPP-002} we obtain that $\sum_{i=1}^{k}\scal{u_{i}}{v_{i}}\le \sum_{j=1}^{k}\cos \fra_{j}$ while \eref{eq-DPP-003} becomes 
\begin{align*}
\frac{1}{2}\sum_{i=1}^{k}\norm{u_{i}-v_{i}}^{2}\ge \sum_{i=1}^{k}\cro{1-\cos \fra_{i}-\frac{1}{2}\pa{\norm{u_{i}}^{2}-\norm{v_{i}}^{2}}}\ge \frd^{2}(U,V)-\frac{1}{2}\sum_{i=1}^{k}\norm{u_{i}-v_{i}}^{2},
\end{align*}
by Pythagoras' theorem. This proves \eref{eq-DPP-aff01}. 

\subsection{Proof of Proposition~\ref{prop-h-DPP}}\label{pf-prop-h-DPP}
For $\gr\in \gS$, let $p(\gr)$ denote the density on $\cJ$  (with respect to the counting measure $\frm$) defined by $J\mapsto p_{J}(\gr)=\prod_{j\in J}\sin^{2}(r_{j})\prod_{j\not\in J}\cos^{2}(r_{j})$.
We observe that the Hellinger affinity between $p(\gs)\cdot\frm$ and $p(\gt)\cdot\frm$ is 
\begin{align*}
\rho(p(\gs)\cdot\frm,p(\gt)\cdot\frm)
&=\sum_{J\in\cJ}\prod_{j\in J}\sin(s_{j})\sin(t_{j})\prod_{j\not \in J}\cos(s_{j})\cos(t_{j})\\
&=\prod_{j\in \frJ}\pa{\sin(s_{j})\sin(t_{j})+\cos(s_{j})\cos(t_{j})}=\prod_{j\in \frJ}\cos(s_{j}-t_{j}).
\end{align*}
Since $\cos(s_{j}-t_{j})\in [0,1]$ for every $s_{j},t_{j}\in [0,\pi/2]$ and $j\in\{1,\ldots,k\}$,  the square Hellinger distance between these probabilities is 
\begin{align*}
h^{2}(p(\gs)\cdot\frm,p(\gt)\cdot\frm)&=1-\rho(p(\gs),p(\gt))=1-\prod_{j\in \frJ}\cos(s_{j}-t_{j})\le \sum_{j\in \frJ}\pa{1-\cos(s_{j}-t_{j})}\\
&= 2\sum_{j\in \frJ}\sin^{2}\pa{\frac{s_{j}-t_{j}}{2}}\le \frac{1}{2}\sum_{j\in \frJ}\ab{s_{j}-t_{j}}^{2}.
\end{align*}
Let $\gr\in\{\gs,\gt\}$. It follows from Proposition~\ref{prop-DPP} that 
\begin{align*}
\sum_{J\in \cJ}p_{J}(\gr)h^{2}\pa{P_{U_{J}}P_{V_{J}}}&\le\frac{1}{2}\sum_{J\in \cJ}p_{J}(\gr)\sum_{j\in \frJ}\norm{u_{j}-v_{j}}^{2}=\frac{1}{2}\sum_{j\in \frJ}\norm{u_{j}-v_{j}}^{2}\sum_{J\in\cJ}p_{J}(\gr)\1_{J}(j)\\
&=\frac{1}{2}\sum_{j\in \frJ}\norm{u_{j}-v_{j}}^{2}\cro{\sum_{J\in\cJ}\pa{\prod_{j'\in J}\sin^{2}(r_{j'})\prod_{j'\not \in \frJ }\cos^{2}(r_{j'})}\1_{ J}(j)}\\
&=\frac{1}{2}\sum_{j\in \frJ}\norm{u_{j}-v_{j}}^{2}\cro{\sin^{2}(r_{j})\prod_{j'\in \frJ, j'\ne j}\pa{\sin^{2}(r_{j'})+\cos^{2}(r_{j'})}}\\
&=\frac{1}{2}\sum_{j\in \frJ}\sin^{2}(r_{j})\norm{u_{j}-v_{j}}^{2}.
\end{align*}
Inequalities in \eref{eq-prop-h-DPP0} are consequences of the following proposition. 
\begin{prop}\label{prop-bar-corr}
Let $(T,\cT,\frm)$ be a measure space, $p,q$ two densities on $(T,\cT,\frm)$ and $\gf=(f_{t})_{t\in T}$, $\gg=(g_{t})_{t\in T}$ two families of probability densities on a measure space $(E,\cE,\mu)$. We denote by
\[
M_{\bs{\phi},r}=\cro{\int_{T}\phi_{t}r(t)d\frm(t)}\cdot\mu\quad \text{for $\bs{\phi}\in\{\gf,\gg\}$ and $r\in\{p,q\}$}
\]
the probability associated with the mixture of the densities $\gh=(h_{t})_{t\in T}$ according to the distribution $r\cdot\frm$ on $T$. 
Then, 
\begin{align*}
h^{2}\pa{M_{\gf,p},M_{\gf,q}}&\le h^{2}(p\cdot\frm,q\cdot\frm)\\
h^{2}\pa{M_{\gf,p},M_{\gg,p}}&\le \int_{T}h^{2}\pa{f_{t}\cdot\mu,g_{t}\cdot \mu}p(t)d\frm(t)\\
h^{2}\pa{M_{\gf,p},M_{\gg,q}}&\le 2\cro{h^{2}(p\cdot\frm,q\cdot\frm)+\inf_{r\in\{p,q\}}\int_{T}h^{2}\pa{f_{t}\cdot\mu,g_{t}\cdot \mu}r(t)d\frm(t)}.
\end{align*}
\end{prop}
\begin{proof}[Proof of Proposition~\ref{prop-bar-corr}]
Note that by the Fubini theorem and the Cauchy--Schwarz inequality
\begin{align*}
\rho\pa{p\cdot \frm,q\cdot\frm}&=\int_{T}\sqrt{p(t)q(t)}d\frm(t)=\int_{E}\cro{\int_{T}f_{t}(x)\sqrt{p(t)q(t)}d\frm(t)}d\mu(x)\\
&\le \int_{E}\sqrt{\int_{T}f_{t}(x)p(t)d\frm(t)}\sqrt{\int_{T}f_{t}(x)q(t)d\frm(t)}d\mu(x)=\rho\pa{M_{\gf,p},M_{\gf,q}}
\end{align*}
which leads to the first inequality. For the second one, using again the Fubini theorem and the Cauchy--Schwarz inequality 
\begin{align*}
\int_{T}\rho\pa{f_{t}\cdot\mu,g_{t}\cdot \mu}p(t)d\frm(t)&=\int_{T}\cro{\int_{E}\sqrt{f_{t}(x)g_{t}(x)}d\mu(x)}p(t)d\frm(t)\\
&\le \int_{E}\cro{\int_{T}\sqrt{f_{t}(x)p(t)}\sqrt{g_{t}(x)p(t)}d\frm(t)}d\mu(x)\\
&\le \rho\pa{M_{\gf,p},M_{\gg,p}},
\end{align*}
which yields the result. The last inequality is a consequence of the triangle inequality and the fact that $(a+b)^{2}\le 2a^{2}+2b^{2}$ for $a,b\ge 0$.
\end{proof}

Inequality \eref{eq-prop-h-DPP01b} is a direct consequence of the first inequality in \eref{eq-prop-h-DPP0} (completing the element $u_{j}$ for $j\in J$ into an orthonormal system of $\H$ if ever necessary). 

Let us now turn to the proof of \eref{eq-prop-h-DPP02b} taking for the sake of simplicity $J=\{1,\ldots,k\}$. For $j\in\{1,\ldots,k\}$, let $v_{j}$ be the orthogonal projection of $u_{j}$ onto $E$. We may proceed as in the proof of Proposition~\ref{prop-GS}. We use the same notation as in that proof. The first part of the proof does not use the fact that the $v_{j}$ are unit vectors. We write $v_{j}=\sum_{i=1}^{k}\sigma_{i}P_{i,j}w_{j}$ and $\overline v_{j}=\sum_{i=1}^{k}P_{i,j}w_{j}$ for $j\in\{1,\ldots,k\}$. The vectors $\overline v_{j}$ are orthonormal, hence $\overline \gv=(\overline v_{1},\ldots,\overline v_{k})\in\U_{k}(E)$, and they satisfy \eref{eq-sigma}. For every $i\in\{1,\ldots,k\}$, $\scal{v_{j}}{w_{i}}=\sigma_{i}P_{i,j}$ and since $P$ is orthogonal
\[
\sigma_{i}^{2}=\sum_{j=1}^{k}\scal{v_{j}}{w_{i}}^{2}=\sum_{j=1}^{k}\scal{u_{j}}{w_{i}}^{2}\le \norm{w_{i}}^{2}=1.
\]
Hence, $\sigma_{i}\in [0,1]$ for every $i\in\{1,\ldots,k\}$ and 
\[
\sum_{i=1}^{k}\sigma_{i}^{2}=\sum_{j=1}^{k}\cro{\sum_{i=1}^{k}\scal{v_{j}}{w_{i}}^{2}}=\sum_{j=1}^{k}|v_{j}|^{2}.
\]
We deduce that 
\begin{align*}
\sum_{j=1}^{k}\norm{u_{j}-\overline v_{j}}^{2}&=2\pa{k-\sum_{j=1}^{k}\scal{u_{j}}{\overline v_{j}}}=2\pa{k-\sum_{j=1}^{k}\scal{v_{j}}{\overline v_{j}}}\\
&=2\pa{k-\sum_{j=1}^{k}\sigma_{j}}\le 2\pa{k-\sum_{j=1}^{k}\sigma_{j}^{2}}=2\sum_{j=1}^{k}\pa{\|u_{j}\|^{2}-\|v_{j}\|^{2}}\\
&=2\sum_{j=1}^{k}\|u_{j}-v_{j}\|^{2}=2\sum_{j=1}^{k}\inf_{v\in E}\|u_{j}-v\|^{2}.
\end{align*}
Finally, we obtain \eref{eq-prop-h-DPP02b} with $\overline \gu=\overline\gv$ by applying this bound to the right-hand side of the second inequality of \eref{eq-prop-h-DPP0} noticing that $\sin^{2} t_{j}$ is bounded by 1 for $j\in J$ and vanishes for $j\not \in J$. 

\subsection{Proof of Proposition~\ref{prop-produit}}\label{pf-prop-produit}
Inequality \eref{eq-prop-produit} is a straightforward consequence of Proposition~\ref{prop-Prop}-\ref{Prop-iii} and \eref{eq-ellmaj}. Let us now prove \eref{eq-prop-produit1}.  With no loss of generality, we may assume that $J=\{1,\ldots,k\}$ with $k\ge 1$ and identify $\gS_{J}$ with $[0,\pi/2]^k$. Using that 
\begin{align*}
B(r)&=\ac{\gt\in \gS_{J}, \sum_{j\in J}\ab{s_{j}-t_{j}}^{2}\le \frac{r}{n}}\supset \ac{\gt\in \gS_{J}, \max_{j\in J}\ab{s_{j}-t_{j}}<\sqrt{\frac{r}{nk}}},
\end{align*}
for $r=nk(\pi z/2)^{2}$ with $z\ge 0$ we obtain that 
\begin{align*}
\gpi_{1}(B(r))&\ge \pa{\frac{2}{\pi}}^{k}\prod_{j=1}^{k}\cro{\pa{s_{j}+\frac{\pi z}{2}}\wedge \frac{\pi}{2}-\pa{s_{j}- \frac{\pi z}{2}}_{+}}\ge \pa{z\wedge 1}^{k}.
\end{align*}
For $z\ge 1$, $\gpi_{1}(B(r))=1$ and the inequality $-\log \gpi_{1}(B(r))\le \gamma r$ is therefore satisfied. For $z<1$, this inequality is satisfied if $az^{2}+b\log z^{2}+c\ge 0$ with $a=\gamma nk\pi^{2}/4$, $b=k/2$ and $c=0$. By applying Lemma~\ref{res-ine} we obtain that $-\log \gpi_{1}(B(r))\le \gamma r$ as soon as $z^{2}\ge (\gamma n\pi^{2}/2)^{-1}\log_{+}(\gamma n\pi^{2}/2)$ which leads to \eref{eq-prop-produit1}.

Let us now turn to the proof of \eref{eq-prop-produit2}. We denote by $\S_{E}$ the unit sphere of $E$ and write for short $\frm_{k}$ for $\frm_{k}^{E}$ for every $k\ge 1$. In particular, $\frm_{1}$ is the uniform distribution on $\S_{E}$. Let us first consider the situation where $D=1$ and $E$ is the linear span of the unit vector $\overline v$. Then $k=1$, $\S_{E}=\{-\overline v,+\overline v\}$, $\U_{1}(E)=\{(-\overline v),(+\overline v)\}$, and for every $r>0$
\begin{align*}
\frm_{1}\pa{\ac{\overline \gv\in\U_{1}(E),\;\sum_{j=1}^{k}\norm{u_{j}-\overline v_{j}}^{2}\le \frac{r}{n}}}\ge \frac{1}{2}.
\end{align*}
The right-hand side is at least $e^{-\gamma r}$ for $r\ge (\log 2)/\gamma$ hence for $r\ge R_{2}(\gpi_{2},\gu)$. 

Let us now assume that $D\ge 2$. We deduce from Proposition~\ref{prop-GS} that for every $\gv\in \S_{E}^{k}$, there exists $\overline \gv\in\U_{k}(E)$ such that 
\begin{align*}
\sum_{j=1}^{k}\norm{u_{j}-\overline v_{j}}^{2}\le 2\sum_{j=1}^{k}\norm{u_{j}-v_{j}}^{2}+2\sum_{j=1}^{k}\norm{v_{j}-\overline v_{j}}^{2}\le 4\sum_{j=1}^{k}\norm{u_{j}-v_{j}}^{2},
\end{align*}
since $\gu\in\U_{k}(E)$. Since $\overline v$ has distribution $\frm_{k}$ when $v_{1},\ldots,v_{k}$ are i.i.d.\ and uniformly distributed on the sphere, we derive that 
\begin{align*}
\frm_{k}\pa{\ac{\overline \gv\in\U_{k}(E),\;\sum_{j=1}^{k}\norm{u_{j}-\overline v_{j}}^{2}\le \frac{r}{n}}}&\ge \frm_{1}^{\otimes k}\pa{\ac{\gv\in \S_{E}^{k},\;\sum_{j=1}^{k}\norm{u_{j}-v_{j}}^{2}\le \frac{r}{4n}}}\\
&\ge \prod_{j=1}^{k}\frm_{1}\pa{\ac{v\in\S_{E},\;\norm{u_{j}-v}\le \sqrt{\frac{r}{4nk}}}}.
\end{align*}
Applying Lemma~7 in Baraud and Romon~\citeyearpar{Baraud:2026aa},  we know that for all $t\in [0,1]$, 
\begin{align*}
\frm_{1}\pa{\ac{v\in\S_{E},\;\norm{u_{j}-v}\le t}}\ge \frac{2^{-D/2}}{\sqrt{\pi D}}t^{D-1}\ge \pa{\frac{t}{\sqrt{2\pi}}}^{D}.
\end{align*}
This implies that for $0< r\le 4nk$, 
\begin{align}\label{eq-UE-min}
-\log \frm_{k}\pa{\ac{\overline \gv\in\U_{k}(E),\;\sum_{j=1}^{k}\norm{u_{j}-\overline v_{j}}^{2}\le \frac{r}{n}}}\le \frac{kD}{2}\log\pa{\frac{8\pi nk}{r}}
\end{align}
and for $r>4nk$, 
\begin{align}\label{eq-UE-min2}
-\log \frm_{k}\pa{\ac{\overline \gv\in\U_{k}(E),\;\sum_{j=1}^{k}\norm{u_{j}-\overline v_{j}}^{2}\le \frac{r}{n}}}\le \frac{kD}{2}\log\pa{2\pi}\le \frac{kD}{2}\log\pa{4\pi}.
\end{align}
We deduce from Lemma~\ref{res-ine} (applied to $x=r$, $a=\gamma$, $b=kD/2$ and $c=-b\log(8\pi nk)$) that the right-hand side of \eref{eq-UE-min} is not larger than $\gamma r$ for  
\[
r\ge \frac{kD}{2\gamma}\log_{+}\pa{\frac{16\pi \gamma n}{D}}\quad \text{while for}\quad r\ge \frac{kD}{2\gamma}\log(4\pi)
\]
that of \eref{eq-UE-min2} does not exceed $\gamma r$. This finishes the proof of \eref{eq-prop-produit2}.

\subsection{Proof of Proposition~\ref{prop-approx}}\label{pf-prop-approx}
Let $\mu$ be the Lebesgue measure on the torus $\T^{d}$. Note that for every $\gk\in\Z_{0}^{d}$ and $\delta>0$
\begin{align*}
\frac{1}{(2\delta)^{d}}\int_{[-\delta,\delta]^{d}}\cos\pa{2\pi\scal{\gk}{x}}d\mu(x)&=\frac{1}{(2\delta)^{d}}\Re\cro{\int_{[-\delta,\delta]^{d}}e^{2i\pi\scal{\gk}{x}}d\mu(x)}\\
&=\frac{1}{(2\delta)^{d}}\Re\cro{\prod_{j=1}^{d}\int_{-\delta}^{\delta}e^{2i\pi k_{j}x_{j}}dx_{j}}=\prod_{j=1}^{d}\frac{\sin(2\pi k_{j}\delta)}{2\pi k_{j}\delta},
\end{align*}
with the convention $(\sin 0)/0=1$. For every $x\in [-\delta,\delta]^d$, 
\begin{align*}
w(\delta)\ge \overline K(0)-\overline K(x)&=2\sum_{\gk\in \Z_{0}^{d}}\overline K_{\gk}\cro{1-\cos\pa{2\pi\scal{\gk}{x}}}\ge 2\sum_{\gk\not \in I}\overline K_{\gk}\cro{1-\cos\pa{2\pi\scal{\gk}{x}}}.
\end{align*}
We deduce that 
\begin{align*}
w(\delta)&\ge \overline K(0)-\frac{1}{(2\delta)^{d}}\int_{[-\delta,\delta]^{d}}\overline K(x)d\mu\ge 2\sum_{\gk\not \in I}\overline K_{\gk}\cro{1-\prod_{j=1}^{d}\ab{\frac{\sin(2\pi k_{j}\delta)}{2\pi k_{j}\delta}}}.
\end{align*}
Note that for $u>0$, $|\sin u|/u\le 1$. If $\gk\not\in I$, there exists $j'\in\{1,\ldots,d\}$ such that $|k_{j'}|>D$, which implies that 
\[
\prod_{j=1}^{d}\ab{\frac{\sin(2\pi k_{j}\delta)}{2\pi k_{j}\delta}}\le \frac{1}{2\pi D\delta}\quad \text{hence}\quad w(\delta)\ge \cro{\overline K(0)-\overline K_{I}(0)}\pa{1-\frac{1}{2\pi D\delta}}.
\]
This yields \eref{eq-approx}. When $w(z)\le Lz^{\alpha}$, we obtain that the minimum of the mapping $\delta\mapsto \pa{1-(2\pi D\delta)^{-1}}^{-1}(L\delta^{\alpha})$ is achieved for $\delta=(1+\alpha)(2\pi D\alpha)^{-1}$ which leads to \eref{eq-approx01}.

\subsection{Proof of Proposition~\ref{prop-nonadapt}}\label{pf-prop-nonadapt}
Let $E=E_{D}$ be the linear span of $1$ and $\frc_{\gk},\frs_{\gk}$ for $\gk\in\sC(D)$. Then $\gv=\fru_{J(D)}=(1,(\frc_{\gk},\frs_{\gk})_{\gk\in\sC(D)})\in\U_{J(D)}(E_{D})$. 
We may equip $\U_{J(D)}(E_{D})$ with the prior $\gpi_{2}$ given by the Dirac mass $\delta_{\fru_{J(D)}}$ at $\fru_{J(D)}$. It is then clear that \eref{def-produit2} is satisfied for all $r\ge 0$ and we may therefore take $R_{2}(\gpi_{2},\fru_{J(D)})=0$.
Let us denote by $\gpi_{I(D)}$ the uniform distribution on $\gS_{I(D)}=[0,\pi/2]\times [0,\pi/2]^{\sC(D)}$. Using that for every $\gt$ in the support of $\gpi_{1,J(D)}$ in $\gS_{J(D)}$,
\begin{align*}
\sum_{j\in J(D)}\ab{s_{j}-t_{j}}^{2}&=\ab{s_{0}-t_{0}}^{2}+\sum_{\gk\in\sC(D)}\sum_{\epsilon\in\{0,1\}}\ab{s_{\gk,\epsilon}-t_{\gk,\epsilon}}^{2}\\
&=\ab{s_{0}-t_{0}}^{2}+2\sum_{\gk\in\sC(D)}\ab{s_{\gk,0}-t_{\gk,0}}^{2}\\
&\le 2\ab{s_{0}-t_{0}}^{2}+2\sum_{\gk\in\sC(D)}\ab{s_{\gk,0}-t_{\gk,0}}^{2}
\end{align*}
we derive from the definition of $\gpi_{1,J(D)}$ that for every $r\ge 0$,
\begin{align*}
&\gpi_{1,J(D)}\pa{\ac{\gt\in\gS_{J(D)},\;\sum_{j\in J(D)}\ab{s_{j}-t_{j}}^{2}\le \frac{r}{n}}}\\
&\ge \gpi_{I(D)}\pa{\ac{\gt\in\gS_{I(D)},\; \ab{s_{\0}-t_{\0}}^{2}+\sum_{\gk\in\sC(D)}\ab{s_{\gk,0}-t_{\gk}}^{2}\le \frac{r}{2n}}}.
\end{align*}
Since $e^{-\gamma r/2}\ge e^{-\gamma r}$, we deduce that $\gpi_{1,J(D)}$ satisfies \eref{def-produit1} for $r\ge 2 R_{1}(\gpi_{I(D)},(s_{0},(s_{\gk,0})_{\gk\in\sC(D)}))$,  which, with \eref{eq-prop-produit1}, implies that we may take 
\[
R_{1}(\gpi_{1,J(D)},\gs_{J(D)})= \frac{1+|\sC(D)|}{\gamma}\log_{+}\pa{\frac{\gamma n\pi^{2}}{2}}\le \frac{(1+2D)^{d}}{\gamma}\log_{+}\pa{\frac{\gamma n\pi^{2}}{2}}.
\]
We conclude by applying Corollary~\ref{cor-su}.

\subsection{Proof of Proposition~\ref{prop-minh}}\label{pf-prop-minh}

As a consequence of \eref{eq-KFour}, let us first note that
\[
\ab{K(x,y)}\le \overline K_{\0}+2\sum_{\gk\in \Z_{0}^{d}}\overline  K_{\gk}\le \overline K(0)\quad \text{for every $x,y\in \T^{d}$.}
\]
Besides, since $K$ is positive semi-definite with diagonal entries $\overline K(0)$, we deduce from Hadamard's inequality that for every $\ell \ge 1$ and $(x_{1},\ldots,x_{\ell})\in(\T^{d})^{\ell}$
\begin{equation}\label{eq-majrhol}
0\le \varrho_{\ell}(x_{1},\ldots,x_{\ell})\le\overline K^{\ell}(0).
\end{equation}

Let  $X,Y$ be two determinantal processes with distribution $P_{\overline K}$ and $P_{\overline L}$ respectively. 
We set  
$\overline \Delta(0)=\ab{\overline K(0)-\overline L(0)}$ and given a finite point process $Z$ on $\T^{d}$,
\[
T_{0}(Z)=Z(\T^{d}) \quad \text{and}\quad T_{1}(Z)=\int_{\T^{d}\times \T^{d}}F(x-y)dZ(x)dZ(y)
\]
with $F(u)=[\overline K^{2}(u)-\overline L^{2}(u)]\1_{u\ne 0}$ for $u\in\T^{d}$. In particular, $\norm{F}_{2}=\|\overline K^2-\overline L^2\|_{2}$. 

By \eref{eq-PDP-f0},   
\begin{align*}
\E\cro{T_{0}(X)}=\int_{\T^{d}}K(x,x)d\mu(x)=\overline K(0)
\end{align*}
and by using \eref{eq-PDP-f}, 
\begin{align*}
\E\cro{T_{0}^{2}(X)}&=\E\cro{\int_{\T^{d}\times \T^{d}}\1_{x\ne y}dX(x)dX(y)+\int_{\T^{d}}dX(x)}\\
&=\int_{\T^{d}\times \T^{d}}\cro{\overline K^{2}(0)-\overline K^{2}(x-y)}d\mu(x)d\mu(y)+\overline K(0)\\
&=\overline K^{2}(0)+\overline K(0)-\int_{\T^{d}}\overline K^{2}(u)d\mu(u)=\overline K^{2}(0)+\overline K(0)-\|\overline K\|_{2}^{2}.
\end{align*}
Consequently, $\Var[X(\T^{d})]= \overline K(0)-\|\overline K\|_{2}^{2}$. Similarly, we obtain that $\Var[Y(\T^{d})]= \overline L(0)-\|\overline L\|_{2}^{2}$ and provided that $\overline \Delta(0)>0$, we deduce from \eref{eq-lem-h0} that
\begin{align}
h^{2}\pa{P_{\overline K},P_{\overline L}}&\ge \frac{\pa{\overline K(0)-\overline L(0)}^{2}}{4\pa{\overline K(0)-\|\overline K\|_{2}^{2}}+4\pa{\overline L(0)-\|\overline L\|_{2}^{2}}+2\pa{\overline K(0)-\overline L(0)}^{2}}\nonumber\\
&\ge \frac{\overline \Delta^{2}(0)}{8 M(0)+2\overline \Delta^{2}(0)}\ge \frac{1}{10}\cro{\frac{\overline \Delta^{2}(0)}{M(0)}\wedge 1}.\label{eq-minH}
\end{align}

Let us now note that 
\begin{align*}
\E\cro{T_{1}(X)}&=\int_{\T^{d}\times \T^{d}}F(x-y)\cro{\overline K^{2}(0)-\overline K^{2}(x-y)}d\mu(x)d\mu(y)\\
&=\overline K^{2}(0)\int_{\T^{d}}F(u)d\mu(u)-\int_{\T^{d}}F(u)\overline K^{2}(u)d\mu(u).
\end{align*}
The same equality holds for $\E\cro{T_{1}(Y)}$ with $\overline L$ in place of $\overline K$. Hence, with our definition of $F$, we obtain that 
\begin{align*}
&\E\cro{T_{1}(X)}-\E\cro{T_{1}(Y)}\\
&=\int_{\T^{d}}F(u)d\mu(u)\pa{\overline K^{2}(0)-\overline L^{2}(0)}+\int_{\T^{d}}F(u)\pa{\overline L^{2}(u)-\overline K^{2}(u)}d\mu(u)\\
&=\cro{\int_{\T^{d}}F(u)d\mu(u)}\pa{\overline K^{2}(0)-\overline L^{2}(0)}-\norm{F}_{2}^{2}
\end{align*}
In particular, 
\begin{equation}\label{eq-minDelta}
\ab{\E\cro{T_{1}(X)}-\E\cro{T_{1}(Y)}}\ge \norm{F}_{2}^{2}-\overline \Delta(0)\cro{\overline K(0)+\overline L(0)}\norm{F}_{2}.
\end{equation}

Let us now turn to the calculation of $\E\cro{T_{1}^{2}(X)}$. Since $F(0)=0$ and $F$ is even, the following decomposition holds 
\begin{align*}
T_{1}^{2}(X)&=\pa{\int_{\T^{d}\times\T^{d}}F(x-y)dX(x)dX(y)}^{2}\\
&=\int_{(\T^{d})^4}F(x_{1}-x_{2})F(x_{3}-x_{4})\prod_{i=1}^{4}dX(x_{i})\\
&=\int_{(\T^{d})^4}F(x_{1}-x_{2})F(x_{3}-x_{4})\1_{|\{x_{1},x_{2},x_{3},x_{4}\}|=4}\prod_{i=1}^{4}dX(x_{i})\\
&\quad +4\int_{(\T^{d})^3}F(x_{1}-x_{2})F(x_{1}-x_{3})\1_{|\{x_{1},x_{2},x_{3}\}|=3}\prod_{i=1}^{3}dX(x_{i})\\
&\quad +2\int_{(\T^{d})^2}F^{2}(x_{1}-x_{2})\prod_{i=1}^{2}dX(x_{i}).
\end{align*}
It follows thus from \eref{eq-PDP-f} that 
\begin{align*}
\E\cro{T_{1}^{2}(X)}&=\int_{(\T^{d})^4}F(x_{1}-x_{2})F(x_{3}-x_{4})\varrho_{4}(x_{1},x_{2},x_{3},x_{4})\prod_{i=1}^{4}d\mu(x_{i})\\
&+4\int_{(\T^{d})^3}F(x_{1}-x_{2})F(x_{1}-x_{3})\varrho_{3}(x_{1},x_{2},x_{3})\prod_{i=1}^{3}d\mu(x_{i})\\
&+2\int_{(\T^{d})^2}F^{2}(x_{1}-x_{2})\varrho_{2}(x_{1},x_{2})\prod_{i=1}^{2}d\mu(x_{i})
\end{align*}
Using \eref{eq-majrhol} and the following equalities and inequalities 
\begin{align*}
\int_{(\T^{d})^4}\ab{F(x_{1}-x_{2})F(x_{3}-x_{4})}\prod_{i=1}^{4}d\mu(x_{i})&=\pa{\int_{\T^{d}}|F(u)|d\mu(u)}^{2}\le\|F\|_{2}^{2}\\
\int_{(\T^{d})^3}\ab{F(x_{1}-x_{2})F(x_{1}-x_{3})}\prod_{i=1}^{3}d\mu(x_{i})&=\pa{\int_{\T^{d}}|F(u)|d\mu(u)}^{2}\le\|F\|_{2}^{2}\\
\int_{(\T^{d})^2}F^{2}(x_{1}-x_{2})\prod_{i=1}^{2}d\mu(x_{i})&=\|F\|_{2}^{2}
\end{align*}
we obtain that 
\begin{align}
\E\cro{T_{1}^{2}(X)}&\le \pa{\overline K^{4}(0)+4\overline K^{3}(0)+2\overline K^{2}(0)}\|F\|_{2}^{2}\le 7\overline K^{2}(0)\pa{\overline K^{2}(0)\vee 1}\|F\|_{2}^{2}.\label{eq-minT2}
\end{align}
The same inequality is satisfied for $\E\cro{T_{1}^{2}(Y)}$ with $\overline L$ in place of $\overline K$.

Let us assume $M(0)>0$ and set 
\[
c=\frac{\sqrt{40M(0)}}{\sqrt{56}(\overline M(0)\vee 1)+\sqrt{40M(0)}}\in (0,1)
\]
 If $\overline \Delta(0)\le c\cro{\overline K(0)+\overline L(0)}^{-1}\norm{F}_{2}$, we deduce from \eref{eq-minDelta} that  
\begin{align*}
\ab{\E\cro{T_{1}(X)}-\E\cro{T_{1}(Y)}}\ge (1-c)\norm{F}_{2}^{2}.
\end{align*}
Provided that  $\norm{F}_{2}^{2}>0$, we deduce from \eref{eq-minT2} and \eref{eq-lem-h1} that in this case
\begin{align}
h^{2}(P_{\overline K},P_{\overline L})&\ge \frac{(1-c)^{2}\norm{F}_{2}^{2}}{28\cro{\overline K^{2}(0)\pa{\overline K^{2}(0)\vee 1}+\overline L^{2}(0)\pa{\overline L^{2}(0)\vee 1}}}\ge \frac{(1-c)^{2}\norm{F}_{2}^{2}}{56\overline M^{2}(0)\cro{\overline M^{2}(0)\vee 1}}\label{eq-0001}\\
&= \frac{\norm{F}_{2}^{2}}{\overline M^{2}(0)\pa{\sqrt{56}(\overline M(0)\vee 1)+\sqrt{40M(0)}}^{2}}\ge  \frac{\norm{F}_{2}^{2}}{191\overline M^{2}(0)\pa{\overline M^{2}(0)\vee 1}}\nonumber
\end{align}
This inequality still holds when $\norm{F}_{2}^{2}=0$. 

If now $\overline \Delta(0)> c\pa{\overline K(0)+\overline L(0)}^{-1}\norm{F}_{2}\ge c\|F\|_{2}/[2\overline M(0)]$, we derive from \eref{eq-minH} that 
\begin{align*}
h^{2}(P_{\overline K},P_{\overline L})&\ge  \frac{\overline \Delta^{2}(0)}{10M(0)}\wedge \frac{1}{10}\ge \frac{c^{2}\|F\|_{2}^{2}}{40\overline M^{2}(0)M(0)}\wedge \frac{1}{10}=\frac{(1-c)^{2}\norm{F}_{2}^{2}}{56\overline M^{2}(0)\cro{\overline M^{2}(0)\vee 1}}\wedge \frac{1}{10}.
\end{align*}
with our choice of $c$. We conclude that 
\begin{align*}
h^{2}(P_{\overline K},P_{\overline L})\ge \pa{ \frac{\overline \Delta^{2}(0)}{10M(0)}\vee \frac{\norm{F}_{2}^{2}}{191\overline M^{2}(0)\pa{\overline M^{2}(0)\vee 1}}}\wedge \frac{1}{10}
\end{align*}
which is \eref{eq-prop-minh} when $M(0)>0$. When $M(0)=0$, $\Var(X(\T^{d}))= \Var(Y(\T^{d}))=0$ and we deduce from \eref{eq-lem-h0} that when $\E\cro{X(\T^{d})}\ne \E\cro{Y(\T^{d})}$, $h^{2}(P_{\overline K},P_{\overline L})\ge 1/2$. If $\E\cro{X(\T^{d})}= \E\cro{Y(\T^{d})}$, $\overline \Delta(0)=0$ and \eref{eq-0001} applies with $c=0$. This proves \eref{eq-prop-minh01}.

\subsection{Proof of Theorem~\ref{thm-borne}}\label{pf-thm-borne}
The function $\overline K$ being continuous at 0, $w(z)\to 0$ as $z\to 0$ and $B_{n}$ thus converges to 0 as $n$ becomes large. Setting 
$\delta_{n}=10C(B_{n}+\xi/n)< 1\wedge \overline K(0)$,  we know from Proposition~\ref{prop-adapt} that on a set $\Omega'$ whose probability is not smaller than $1-2e^{-\xi}$, 
\begin{equation}\label{eref-deltan}
h^{2}(P_{\overline K},P_{\widehat K})\le \frac{\delta_{n}}{10}< \frac{1}{10}.
\end{equation}
Let $M(0)$ and $\overline M(0)\ge M(0)$ as defined in Proposition~\ref{prop-minh} with $L=\widehat K$. Let us first assume that $M(0)>0$. It then follows from \eref{eq-prop-minh} and \eref{eref-deltan} that on $\Omega'$
\begin{align*}
\frac{|\overline K(0)-\widehat K(0)|^{2}}{\overline K(0)\vee \widehat K(0)}\le \delta_{n}.
\end{align*}
If additionally $\widehat K(0)\ge \overline K(0)$, we get 
\begin{align*}
\widehat K(0)\le \overline K(0)+\sqrt{\widehat K(0)\delta_{n}}\le  \overline K(0)+\frac{1}{2}\widehat K(0)+\frac{1}{2}\delta_{n}\quad \text{that is}\quad \widehat K(0)\le 2\overline K(0)+\delta_{n}
\end{align*}
while $\widehat K(0)<\overline K(0)$ otherwise. We deduce that on $\Omega'$, $\widehat K(0)\le 2\cro{\overline K(0)+\delta_{n}}\le 4\overline K(0)$ when $M(0)>0$. When $M(0)=0$, we necessarily have $\widehat K(0)=\overline K(0)$ by Proposition~\ref{prop-minh}
and \eref{eref-deltan}. The inequality $\widehat K(0)\le 4\overline K(0)$ is therefore still satisfied when $M(0)=0$. We infer from  \eref{eq-prop-minh},\eref{eq-prop-minh01} and \eref{eref-deltan} that on $\Omega'$, 
\[
\norm{\overline K^{2}-\widehat K^{2}}_{2}^{2}\le C'\pa{1\vee \overline K(0)}^{2}\overline K^{2}(0)\delta_{n}
\]
for some positive constant $C'$, which yields \eref{eq-thm-borne}.

\bibliographystyle{apalike}

\end{document}